\documentclass[11pt]{article}

\usepackage{bm}
\usepackage{amsfonts}
\usepackage{amssymb,amsmath,amsthm, amsfonts}
\usepackage{mathrsfs}

\def\sqr#1#2{{\vcenter{\vbox{\hrule height.#2pt
              \hbox{\vrule width.#2pt height#1pt \kern#1pt \vrule width.#2pt}
          \hrule height.#2pt}}}}
\def\sqr#1#2{{\vcenter{\vbox{\hrule height.#2pt
              \hbox{\vrule width.#2pt height#1pt \kern#1pt \vrule width.#2pt}
              \hrule height.#2pt}}}}
\def\3n{\negthinspace \negthinspace \negthinspace }
\def\2n{\negthinspace \negthinspace }
\def\1n{\negthinspace }

\def\={\buildrel \triangle \over =}

\def\limsup{\mathop{\overline{\rm lim}}}

\def\esssup{\mathop{\rm esssup}}

\def\max{\mathop{\rm max}}

\def\exp{\mathop{\rm exp}}
\def\sup{\mathop{\rm sup}}
\def\inf{\mathop{\rm inf}}

\def\({\Big (}
\def\){\Big )}
\def\[{\Big[}
\def\]{\Big]}
\def\be{\begin{equation}}

\def\ee{\end{equation}}

\def\square#1{\vbox{\hrule\hbox{\vrule height#1%
     \kern#1\vrule}\hrule}}
\def\rectangle#1#2{\vbox{\hrule\hbox{\vrule height#1%
     \kern#2\vrule}\hrule}}

\font\tenbb=msbm10 \font\sevenbb=msbm7 \font\fivebb=msbm5

\newfam\bbfam
\scriptscriptfont\bbfam=\fivebb \textfont\bbfam=\tenbb
\scriptfont\bbfam=\sevenbb

\theoremstyle{definition}
\newtheorem{lemma}{Lemma}[section]
\newtheorem{remark}{Remark}[section]

\newtheorem{theorem}{Theorem}[section]

\newtheorem{definition}{Definition}[section]
\newtheorem{proposition}{Proposition}[section]

\makeatletter
   
   \@addtoreset{equation}{section}
\makeatother

\begin{document}

\title{General mean-field BSDEs with diagonally quadratic generators in multi-dimension\thanks{Juan Li is supported by the NSF of P.R. China (NOs. 12031009, 11871037), National Key R and D Program of China (NO. 2018YFA0703900), NSFC-RS (No. 11661130148; NA150344). Qingmeng Wei is supported by the Natural Science Foundation of Jilin Province for Outstanding Young Talents (No. 20230101365JC) and National Natural Science Foundation of China (11971099, 12371443).}}
\author{Weimin Jiang$^{1}$,\,\, Juan Li$^{1,2,}$\footnote{Corresponding authors.},\,\, Qingmeng Wei$^{3, \dag}$ \\
{$^1$\small School of Mathematics and Statistics, Shandong University, Weihai, Weihai 264209, P.~R.~China.}\\
{$^2$\small Research Center for Mathematics and Interdisciplinary Sciences, Shandong University,}\\
{\small Qingdao 266237, P.~R.~China.}\\
{$^3$\small School of Mathematics and Statistics, Northeast Normal University, Changchun 130024, P.~R.~China.}\\
{\footnotesize{\it E-mails: weiminjiang@mail.sdu.edu.cn,\,\ juanli@sdu.edu.cn,\,\ weiqm100@nenu.edu.cn.}}
}
\date{\today}
\maketitle

\textbf{Abstract}. The purpose of this paper is to investigate general mean-field backward stochastic differential equations (MFBSDEs) in multi-dimension with diagonally quadratic generators $f(\omega,t,y,z,\mu)$, that is, the coefficients depend not only on the solution processes $(Y,Z)$, but also on their law $\mathbb{P}_{(Y,Z)}$, as well as have a diagonally quadratic growth in $Z$ and super-linear growth (or even a quadratic growth) in the law of $Z$ which is totally new. We start by establishing through a fixed point theorem the existence and the uniqueness of local solutions in the ``Markovian case'' $f(t,Y_{t},Z_{t},\mathbb{P}_{(Y_{t},Z_{t})})$ when the terminal value is bounded. Afterwards, global solutions are constructed by stitching local solutions. Finally, employing the $\theta$-method, we explore the existence and the uniqueness of global solutions for diagonally quadratic mean-field BSDEs with convex generators, even in the case of unbounded terminal values that have exponential moments of all orders. These results are extended to a Volterra-type case where the coefficients can even be of quadratic growth with respect to the law of $Z$.

\textbf{Keywords}. Mean-field BSDEs; Diagonally quadratic generator; BMO martingale; Convex generator; Unbounded terminal value

\section{Introduction}
\hspace{2em}The origin of linear BSDEs can be traced back to their role as adjoint equation in stochastic control problems, which was studied by Bismut \cite{B73} in 1973. Later on, Pardoux and Peng \cite{PP90} were the first to introduce nonlinear BSDEs and to prove the existence and the uniqueness of adapted solutions under the condition that the generator is Lipschitz continuous and the terminal value is square integrable. The form of the equation writes as follows:
\begin{equation}\label{eq 1.1}\tag{1.1}
	Y_{t}=\xi+\int_{t}^{T}f(s,Y_{s},Z_{s})ds-\int_{t}^{T}Z_{s}dW_s,\ t\in[0,T].
\end{equation}

Since their pioneering paper the theory of BSDEs has attracted a lot of researchers because of its large fields of applications, e.g., El Karoui, Peng and Quenez \cite{EP97}, Pardoux and Peng \cite{PP92}, Peng \cite{P92}, etc. Numerous papers have been devoted to improving the existence and the uniqueness result of Pardoux and Peng \cite{PP90} by relaxing the Lipschitz hypothesis on the coefficients. In particular, Lepeltier and San Martin \cite{LM97} considered the existence of a solution for one-dimensional BSDE when the generator is only continuous and has a linear growth.

The study of the solvability of BSDEs with quadratic growth constitutes a significant research field in its own right. These equations arise, by example, in the context of utility optimization problems with exponential utility functions, or alternatively in problems related to risk minimization for the entropic risk measure. Kobylanski \cite{K00} first provided the existence and the uniqueness result for the case of a scalar-valued quadratic BSDE with bounded terminal value. Since then, there have been many further studies on this subject. The investigation carried out by Briand and Hu \cite{BH08} focused on quadratic BSDEs which involve terminal values that are unbounded, as well as generators that are convex. In comparison with the one-dimensional BSDE with quadratic growth the multi-dimensional ones pose a considerable challenge, as certain tools of the one-dimensional case, such as comparison theorem, monotone convergence or Girsanov theorem, are no longer applicable here. Frei and Dos Reis \cite{FD11} constructed a counterexample to illustrate that, even in seemingly benign situations, the existence of global solutions may fail. Nevertheless, we can find in the literature several papers that deal with special cases of multi-dimensional quadratic BSDEs. Tevzadze \cite{T08} studied the existence and the uniqueness results for multi-dimensional quadratic BSDEs (1.1) while assuming that the terminal value has sufficiently small supremum. Cheridito and Nam \cite{CN15} analysed a system of BSDEs with projectable quadratic generators. Hu and Tang \cite{HT16} established the solvability of local and global solutions for systems of BSDEs with diagonally quadratic generators and discussed a nonzero-sum risk-sensitive differential game of stochastic functional differential equations. Xing and \v{Z}itkovi\'{c} \cite{XZ18} obtained a general result in a Markovian setting with weak regularity assumptions on the coefficients. Fan, Hu and Tang \cite{FH23} studied diagonally quadratic BSDEs by weakening the assumptions of Hu and Tang \cite{HT16}.

Mean-field stochastic differential equations, also known as McKean-Vlasov equations, which describe stochastic systems whose evolution is determined by both the microcosmic location and the macrocosmic distribution of the particles, can be dated back to the original work by Kac \cite{K56} on the Boltzmann or the Vlasov equations in the kinetic gas theory in the 1950s. Mean-field game theory started with the groundbreaking contributions of Lasry and Lions \cite{LL07}. Inspired by this, Buckdahn, Djehiche, Li and Peng \cite{BD09} introduced mean-field BSDEs, which can be derived from an approximation involving a weakly interacting $N$-particle system, and investigated it with the help of purely probabilistic methods. Buckdahn, Li and Peng \cite{BL09} studied the existence and the uniqueness of solutions as well as a comparison theorem for mean-field BSDEs when the generator $f$ is Lipschitz continuous and the terminal condition $\xi$ is square integrable. Since then, there has been a growing interest for mean-field SDEs and BSDEs, e.g., Buckdahn, Li, Peng and Rainer \cite{BL17}, Hao, Wen and Xiong \cite{HW22}, Hibon, Hu and Tang \cite{HH22-2}, Li \cite{L18}, Li, Liang and Zhang \cite{LL18}. Recently, Hao, Hu, Tang and Wen \cite{HH22-1} studied the solvability of local and global solutions for one-dimensional mean-field BSDEs with quadratic growth. Tang and Yang \cite{TY23} established the existence and the uniqueness of local and global solutions for multi-dimensional BSDEs whose generators depend on the expectation of both variables.

Motivated strongly by the above works, we study multi-dimensional diagonally quadratic mean-field BSDEs, in which the generator relies on both the solution $(Y,Z)$ and the law of $(Y,Z)$. More precisely, we consider the following multi-dimensional mean-field BSDE
\begin{equation}\label{eq 1.2}\tag{1.2}
	Y_{t}=\xi+\int_{t}^{T}f(s,Y_{s},Z_{s},\mathbb{P}_{(Y_{s},Z_{s})})ds-\int_{t}^{T}Z_{s}dW_s,\ t\in[0,T],
\end{equation}
where the generator $f$ is diagonally quadratic growth with respect to $z$, i.e., the quadratic part of the $i$th component depends only on the $i$th row of $z$, and $\mathbb{P}_{(Y_{s},Z_{s})}$ denotes the joint law of $Y_{s}$ and $Z_{s}$. In this paper, we will focus on the study of the existence and the uniqueness theorem of both local and global solutions for mean-field BSDEs of diagonally quadratic generators. By means of the fixed point argument, we prove the existence and the uniqueness of local solutions. This result will be used for the construction of global solutions. In a next step we solve these equations with unbounded terminal values and convex generators through the use of iterative algorithm, uniform a priori estimates and the $\theta$-method. Finally, we extend the above ``Markovian case'' with driver $f(t,Y_{t},Z_{t},\mathbb{P}_{(Y_{t},Z_{t})})$ to the Volterra-type case $\displaystyle f_{s}(t,Y,Z,\mathbb{P}_{(Y,Z)})=\mathbb{E}\big[g(t,Y_{\cdot\vee t},Z,\mathbb{P}_{(Y_{\cdot\vee t},Z)})|\mathcal{F}_{s}\big]+f(t,Z_{t}),\ s\le t\le T,$ where $f(t,z)$ has a diagonally quadratic growth in $z\in \mathbb{R}^{n\times d}$, and $g$ depends at time $t$ on the future path of $Y$ (on $Y_{\cdot \vee t}=\{Y_{\max(s,t)}\}_{s\in[0,T]}$) and the full path of $Z$ as well as on the joint law of both.

Our first contribution is to generalize non-trivially the existence and the uniqueness result for local solutions in \cite{FH23} to the mean-field case (see, Theorem 3.1), which will be a tool for studying non-zero sum risk-sensitive stochastic differential games in mean-field case later. Our second contribution involves investigating the existence and the uniqueness of global solutions for diagonally quadratic mean-field BSDEs in multi-dimension, which covers the one-dimensional mean-field case in \cite{HH22-1} (see, Theorem 4.1). Our third contribution lies in the generalization of the multi-dimensional diagonally quadratic BSDEs with unbounded terminal values in \cite{FH23} to the mean-field case, where the coefficients also depend on the distribution of $Y$, and we establish the solvability of global solutions (see, Theorem 5.1). Finally, we further prove the existence and the uniqueness of solutions for path-dependent Volterra-type mean-field BSDEs of diagonally quadratic generators with unbounded terminal values, where one of the coefficients can even be of quadratic growth in the law of $Z$ (see, Theorem 5.2). In the existing works for the multi-dimensional case without mean-field term the coefficients can only have a diagonally quadratic growth in $Z$.

The rest of this paper is organized as follows: Section 2 presents notations and recalls some fundamental preliminaries about the properties of BMO martingales. In Section 3, we study the existence and the uniqueness of local solutions  for diagonally quadratic mean-field BSDEs with bounded terminal values by using a fixed point theorem. Section 4 is devoted to the solvability of global solutions for these equations. In Section 5, with the aid of the $\theta$-method, we discuss the existence and the uniqueness results for diagonally quadratic mean-field BSDEs with unbounded terminal values, but possessing exponential moments of all orders. Finally, using the such obtained results, we discuss the Volterra-type mean-field BSDE with a diagonally quadratic growth in $Z$ and quadratic growth in the law of $Z$.
\section{Preliminaries}
\hspace{2em} Throughout this paper, we fix a finite time horizon $T\in (0,+\infty)$ and we let $(\Omega,\mathcal{F},\mathbb{P})$ be a complete probability space endowed with a $d$-dimensional standard Brownian motion $\{W_{t}\}_{t\in[0,T]}$. Let $a\wedge b$ and $a\vee b$ denote the minimum and maximum of two real numbers $a$ and $b$, respectively. Set $a^{+}:=a\vee 0$ and $a^{-}:=-(a\wedge 0)$. Denote by ${\bm{1}}_{A}(\cdot)$ the indicator of set $A$, and sgn$(x):={\bm{1}}_{\{x>0\}}-{\bm{1}}_{\{x\le0\}},\ x\in \mathbb{R}$. For $1\le i\le n$, denote by $z^{i},y^{i}$ and $f^{i}$ respectively the $i$th row of matrix $z\in \mathbb{R}^{n\times d}$, the $i$th component of the vector $y\in \mathbb{R}^{n}$ and the generator $f$. The notation $\delta_{0}$ denotes the Dirac measure at $0$. We suppose that there exists a sub-$\sigma$-field $\mathcal{F}^0\subset \mathcal{F},$  $\mathcal{F}^0$ includes all $\mathbb{P}$-null subsets of $\mathcal{F},$ such that

\indent\quad(i)\ \ the Brownian motion $W$ is independent of $\mathcal{F}^{0};$

\indent\quad(ii)\ \ $\mathcal{F}^0$ is ``rich enough''$,$ i.e.$,$ $\mathcal{P}_2(\mathbb{R}^k)$=$\{\mathbb{P}_\vartheta, \vartheta\in L^2_{\mathcal{F}^{0}}(\Omega;\mathbb{R}^k)\},$ $k\geq 1$.\\
\noindent Here $\mathbb{P}_\vartheta:=\mathbb{P}\circ[\vartheta]^{-1}$ denotes the law of the random variable $\vartheta$ under the probability $\mathbb{P}$. We denote by $\mathbb{F}=(\mathcal{F}_{t})_{t\in[0,T]}$ the filtration generated by $W$ and augmented by $\mathcal{F}^{0}$, i.e., $\mathcal{F}_t:=\sigma\{W_s,\ 0\leq s\leq t\}\vee\mathcal{F}^{0},\ t\in[0,T]$.

The following spaces will be frequently used in this paper.

\noindent $\bullet\ L^{p}_{\mathcal{F}_{t}}(\Omega;\mathbb{R}^{n})\ \text{is the set of}\ \mathbb{R}^{n}\text{-valued},\ \mathcal{F}_{t}\text{-measurable}\ \text{random variables}\ \xi\ \text{such that}\ \Vert\xi\Vert_{L^{p}(\Omega)}:=\big(\mathbb{E}\big[|\xi|^{p}\big]\big)^\frac{1}{p}<+\infty.$

\noindent $\bullet\ L^{\infty}_{\mathcal{F}_{t}}(\Omega;\mathbb{R}^{n})\ \text{is the set of all}\ \mathbb{R}^{n}\text{-valued},\ \mathcal{F}_{t}\text{-measurable}\ \text{random variables}\ \xi\ \text{such that}\ \Vert\xi\Vert_{\infty}:=\esssup\limits_{\omega\in \Omega}|\xi(\omega)|<+\infty.$

\noindent $\bullet\ \mathcal{S}^{p}_{\mathbb{F}}(t,T;\mathbb{R}^{n})\ \text{is the set of}\ \mathbb{R}^{n}\text{-valued}\ \mathbb{F}\text{-adapted}\ \text{continuous processes}\ (\varphi_{s})_{s\in [t,T]}\ \text{with}\ \Vert\varphi\Vert_{\mathcal{S}^{p}_{\mathbb{F}}(t,T)}:=\big(\mathbb{E}\big[\sup\limits_{s\in[t,T]}|\varphi_{s}|^{p}\big]\big)^\frac{1}{p}<+\infty.$

\noindent $\bullet\ \mathcal{S}^{\infty}_{\mathbb{F}}(t,T;\mathbb{R}^{n})\ \text{is the set of}\ \mathbb{R}^{n}\text{-valued}\ \mathbb{F}\text{-adapted}\ \text{continuous processes}\ (\varphi_{s})_{s\in [t,T]}\ \text{with}\ \Vert\varphi\Vert_{\mathcal{S}^{\infty}_{\mathbb{F}}(t,T)}:=\esssup\limits_{(s,\omega)\in [t,T]\times \Omega}|\varphi_{s}(\omega)|<+\infty.$

\noindent $\bullet\ \mathcal{H}^{p}_{\mathbb{F}}(t,T;\mathbb{R}^{n})\ \text{is the set of}\ \mathbb{R}^{n}\text{-valued}\ \mathbb{F}\text{-adapted measurable}\ \text{processes}\ (\psi_{s})_{s\in [t,T]}\ \text{with}\ \Vert\psi\Vert_{\mathcal{H}^{p}_{\mathbb{F}}(t,T)}:=\big(\mathbb{E}\big[\big(\int_{t}^{T}|\psi_{s}|^{2}ds\big)^{\frac{p}{2}}\big]\big)^{\frac{1}{p}}<+\infty.$

\noindent $\bullet\ \mathcal{Z}^{2}_{\mathbb{F}}(t,T;\mathbb{R}^{n})\ \text{is the subset of}\ \mathcal{H}^{2}_{\mathbb{F}}(t,T;\mathbb{R}^{n})\ \text{with}\ \Vert\psi\Vert_{\mathcal{Z}^{2}_{\mathbb{F}}(t,T)}:=\sup\limits_{\tau\in \mathscr{T}[t,T]}\Big\Vert\mathbb{E}_{\tau}\Big[\int_{\tau}^{T}|\psi_{s}|^{2}ds\Big]\Big\Vert_{\infty}^{\frac{1}{2}}<+\infty,$

\noindent where $\mathscr{T}[t,T]$ denotes the set of all $\mathbb{F}$-stopping times $\tau$ with values in $[t,T]$ and $\mathbb{E}_{\tau}$ is the conditional expectation with respect to $\sigma$-field $\mathcal{F}_{\tau}$. We denote by $\mathbb{E}^{\mathbb{Q}}$ and $\mathbb{E}_{\tau}^{\mathbb{Q}}$ the expectation and the conditional expectation, respectively, under the probability measure $\mathbb{Q}$. Moreover, let $M=(M_{t},\mathcal{F}_{t})$ be a uniformly integrable martingale with $M_{0}=0$. For $p\ge 1$, we set
$$\Vert M\Vert_{BMO_{p}(\mathbb{P})}:=\sup\limits_{\tau\in \mathscr{T}[0,T]}\Big\Vert E_{\tau}\left[\big(\langle M \rangle_{T}-\langle M \rangle_{\tau}\big)^{\frac{p}{2}}\right]^{\frac{1}{p}}\Big\Vert_{\infty}< +\infty.$$

\noindent The class $\{M: \Vert M\Vert_{BMO_{p}(\mathbb{P})}<+\infty\}$ is denoted by $BMO_{p}(\mathbb{P})$. Observe that $BMO_{p}(\mathbb{P})$ is a Banach space equipped with the norm $\Vert\cdot\Vert_{BMO_{p}(\mathbb{P})}$. In the sequel for simplicity, we write $BMO(\mathbb{P})$ for the space of $BMO_{2}(\mathbb{P})$. Denote by $\mathscr{E}(M)$ the Dol\'{e}ans-Dade exponential of a one-dimensional local martingale $M$ and by $\mathscr{E}(M)_{t}^{s}$ that of $M_{s}-M_{t}$. We recall that the process $t\rightarrow \int_{0}^{t}Z_{s}dW_{s}$ (denoted by $Z\cdot W$) belongs to $BMO(\mathbb{P})$ if and only if $Z\in \mathcal{Z}_{\mathbb{F}}^{2}(0,T)$. Moreover, in this case,
$$\Vert Z\cdot W\Vert_{BMO(\mathbb{P})}=\Vert Z\Vert_{\mathcal{Z}^{2}_{\mathbb{F}}(0,T)}.$$

\noindent For $p\ge1$, we denote by $\mathcal{P}_{p}(\mathbb{R}^k)$ the family of all probability measures $\mu$ on $(\mathbb{R}^k,\mathcal{B}(\mathbb{R}^k))$ with finite $p$th-order moment, i.e., $\Vert\mu\Vert_{p}:=\big(\int_{\mathbb{R}^{k}}|x|^{p}\mu(dx)\big)^{\frac{1}{p}}<+\infty$, where $\mathcal{B}(\mathbb{R}^k)$ denotes the Borel $\sigma$-field over $\mathbb{R}^{k}$. The space $\mathcal{P}_{p}(\mathbb{R}^k)$ is endowed with the $p$-Wasserstein metric: For $\mu,\nu\in\mathcal{P}_{p}(\mathbb{R}^k),$
\begin{equation*} W_{p}(\mu,\nu):=\inf\Big\{\Big(\int_{\mathbb{R}^k\times\mathbb{R}^k}|x-y|^{p}\rho(dxdy)\Big)^{\frac{1}{p}},\ \rho\in\mathcal{P}_{p}(\mathbb{R}^{k}\times\mathbb{R}^{k})\ \text{with marginals}\ \mu\ \text{and}\ \nu\Big\}.
\end{equation*}


\noindent Finally, we write $Y\in \mathcal{E}(t,T;\mathbb{R}^{n})$ if $\exp\big\{\sup_{s\in[t,T]}|Y_{s}|\big\}\in \mathcal{S}^{p}_{\mathbb{F}}(t,T;\mathbb{R}^{n})$, for all $p\ge1$, and $Z\in \mathcal{M}(t,T;\mathbb{R}^{n\times d})$, if $Z\in \mathcal{H}^{p}_{\mathbb{F}}(t,T;\mathbb{R}^{n\times d})$ for all $p\ge1$.

\begin{definition} A pair of processes $(Y,Z)\in \mathcal{S}^{2}_{\mathbb{F}}(0,T;\mathbb{R}^{n})\times \mathcal{H}^{2}_{\mathbb{F}}(0,T;\mathbb{R}^{n\times d})$ is said to be an adapted solution of mean-field BSDE (\ref{eq 1.2}), if $\mathbb{P}$-almost surely, it satisfies (\ref{eq 1.2}), and a bounded adapted solution if further it belongs to $\mathcal{S}^{\infty}_{\mathbb{F}}(0,T;\mathbb{R}^{n})\times \mathcal{Z}^{2}_{\mathbb{F}}(0,T;\mathbb{R}^{n\times d})$.
\end{definition}

We list below several properties associated with BMO martingales which will be used later. For a more comprehensive analysis, we refer the reader to Kazamaki \cite{K94} and the references therein.

\begin{lemma} (John-Nirenberg Inequality) Let $M\in BMO(\mathbb{P})$ such that $\Vert M\Vert_{BMO(\mathbb{P})}<1$. Then for every stopping time $\tau\in \mathscr{T}[0,T]$, we have
$$\mathbb{E}_{\tau}[\exp\{\langle M\rangle_{T}-\langle M\rangle_{\tau}\}]\le \Big(1-\Vert M\Vert^{2}_{BMO(\mathbb{P})}\Big)^{-1}.$$
\end{lemma}

\begin{lemma} For every $p\in[1,\infty)$, there is a generic constant $L_{p}>0$ such that for all uniformly integrable martingale $M$, the following inequality holds true:
$$\Vert M\Vert_{BMO_{p}(\mathbb{P})}\le L_{p}\Vert M\Vert_{BMO(\mathbb{P})}.$$
\end{lemma}


\begin{lemma} For $K>0$, there exist constants $c_{1}>0$ and $c_{2}>0$ depending only on $K$ such that for every $M\in BMO(\mathbb{P})$, we have for all one-dimensional $N\in BMO(\mathbb{P})$ satisfying $\Vert N \Vert_{BMO(\mathbb{P})}\le K$,
$$c_{1}\Vert M\Vert_{BMO(\mathbb{P})}\le \Vert \widetilde{M}\Vert_{BMO(\widetilde{\mathbb{P}})}\le c_{2}\Vert M\Vert_{BMO(\mathbb{P})},$$
where $\widetilde{M}:=M-\langle M,N\rangle$ and $d\widetilde{\mathbb{P}}:=\mathscr{E}(N)_{0}^{T}d\mathbb{P}$.
\end{lemma}

\section{Local solution with bounded terminal value}

\hspace{2.6em}In this section, our objective is to prove the existence and the uniqueness of local solutions for mean-field BSDE (\ref{eq 1.2}) with diagonally quadratic generators. We rewrite the equation as follows:
\begin{equation}\label{eq 3.1}\tag{3.1}
	Y_{t}=\xi+\int_{t}^{T}f(s,Y_{s},Z_{s},\mathbb{P}_{(Y_{s},Z_{s})})ds-\int_{t}^{T}Z_{s}dW_s,\ t\in[0,T].
\end{equation}

In what follows, we always suppose that $\zeta:[0,T]\times \Omega\rightarrow \mathbb{R}_{+}$ is an $\mathbb{F}$-progressively measurable process, $\psi,\psi_{0}:[0,+\infty)\rightarrow[0,+\infty)$ are two increasing continuous functions, and $\gamma,\gamma_{0},\lambda,K,L$ and $\alpha\in[0,1)$ are positive constants. For the sake of notational simplicity, we shall use $\mu_{1}=\mu(\cdot\times \mathbb{R}^{n\times d})$ and $\mu_{2}=\mu(\mathbb{R}^{n}\times \cdot)$ to denote the marginal distributions of $\mu\in \mathcal{P}_{2}(\mathbb{R}^{n+n\times d})$. Let $\xi:\Omega\rightarrow \mathbb{R}^{n}$ be $\mathcal{F}_{T}$-measurable and $f:\Omega\times[0,T]\times\mathbb{R}^{n}\times\mathbb{R}^{n\times d}\times\mathcal{P}_{2}(\mathbb{R}^{n+n\times d})\rightarrow \mathbb{R}^{n}$ be $\mathbb{F}$-progressively measurable. Moreover, the following assumptions are supposed to be satisfied:\vspace{1em}

\noindent\textbf{(A1)} There exists $M_{1}\ge 0$ such that $\Vert \xi \Vert_{\infty}\le M_{1}$.

\noindent\textbf{(A2)} For $1\le i\le n,$ and $d\mathbb{P}\times dt$-almost all $(\omega,t)\in \Omega\times [0,T],$ there exist constants $\gamma>0,\lambda>0,\gamma_{0}>0,M_{2}>0$ and an $\mathbb{F}$-progressively measurable process $\zeta=(\zeta_{t})$ with $\big\Vert
\int_{0}^{T}\zeta_{t}(\omega)dt\big\Vert_{\infty}\le M_{2}$, such that, for all $y\in \mathbb{R}^{n},$ $z\in \mathbb{R}^{n\times d},$ and $\mu\in \mathcal{P}_{2}(\mathbb{R}^{n+n\times d}),$
\begin{equation*}
|f^{i}(\omega,t,y,z,\mu)|\le \zeta_{t}(\omega)+\psi(|y|)+\frac{\gamma}{2}|z^{i}|^{2}+\lambda \sum\limits_{j\ne i}|z^{j}|^{1+\alpha}+\psi_{0}\big(W_{2}(\mu_{1},\delta_{0})\big)+\gamma_{0}W_{2}(\mu_{2},\delta_{0})^{1+\alpha}.
\end{equation*}

\noindent\textbf{(A3)} For $1\le i\le n,$ and $d\mathbb{P}\times dt$-almost all $(\omega,t)\in \Omega\times [0,T],$ we have that, for all $y,\bar{y}\in \mathbb{R}^{n},$ $z,\bar{z}\in \mathbb{R}^{n\times d},$ and $\mu,\bar{\mu}\in \mathcal{P}_{2}(\mathbb{R}^{n+n\times d}),$
\begin{equation*}
\begin{aligned}
&|f^{i}(\omega,t,y,z,\mu)-f^{i}(\omega,t,\bar{y},\bar{z},\bar{\mu})|\\
&\hspace{-2em}\le\ \psi\big(|y|\vee|\bar{y}|\vee W_{2}(\mu_{1},\delta_{0}) \vee W_{2}(\bar{\mu}_{1},\delta_{0})\big)\\
&\times\bigg[\Big(1+|z|+|\bar{z}|+W_{2}(\mu_{2},\delta_{0})+W_{2}(\bar{\mu}_{2},\delta_{0})\Big)\Big(|y-\bar{y}|+|z^{i}-\bar{z}^{i}|+W_{2}(\mu_{1},\bar{\mu}_{1})\Big)\\
&\hspace{1.5em}+\Big(1+|z|^{\alpha}+|\bar{z}|^{\alpha}+W_{2}(\mu_{2},\delta_{0})^{\alpha}+W_{2}(\bar{\mu}_{2},\delta_{0})^{\alpha}\Big)\Big(\sum\limits_{j\ne i}|z^{j}-\bar{z}^{j}|+W_{2}(\mu_{2},\bar{\mu}_{2})\Big)\bigg].\\[-2mm]
\end{aligned}
\end{equation*}


\begin{remark} The assumptions (A2)-(A3) are more general than the conditions (H1)-(H2) in \cite{FH23} for the non-mean-field case and also than the conditions (H1)-(H2) in \cite{TY23} where the coefficients depend also on the expectation of the solutions. For instance, the generator $f$, specified below, satisfies the assumptions (A2)-(A3), but fails to fulfill the conditions (H1)-(H2) of \cite{TY23}: For example, for $1\le i\le n,\ t\in[0,T],$ $y\in \mathbb{R}^{n},$ $z\in \mathbb{R}^{n\times d},$ and $\mu\in \mathcal{P}_{2}(\mathbb{R}^{n+n\times d}),$ we can consider the function
\begin{equation*}f^{i}(t,y,z,\mu)=(|y|^{2}+\sin|z^{i}|)|z|+|z|^{\frac{4}{3}}+|z^{i}|^{2}+W_{2}(\mu_{1},\delta_{0})^{3}\cos(W_{2}(\mu_{2},\delta_{0}))+W_{2}(\mu_{2},\delta_{0})^{\frac{4}{3}},
\end{equation*}
but \cite{TY23} can not.
\end{remark}

By generalizing the corresponding results of \cite[Proposition 3.2]{HH22-1} for the one-dimensional mean-field case and \cite[Lemma A.1]{FH23} for the multi-dimensional case, we obtain the following crucial extensions for the existence and the uniqueness as well as an a priori estimate for quadratic BSDEs.

\begin{proposition} Given a pair of processes $(U,V)\in \mathcal{S}^{\infty}_{\mathbb{F}}(0,T;\mathbb{R}^{n})\times \mathcal{Z}^{2}_{\mathbb{F}}(0,T;\mathbb{R}^{n\times d})$, assume that  the generator $f:\Omega\times[0,T]\times \mathbb{R}^{1\times d}\rightarrow \mathbb{R}$ and the terminal value $\xi$ satisfy the conditions:

\noindent(i) There exists $M_{1}\ge 0$ such that $\Vert \xi \Vert_{\infty}\le M_{1}$;

\noindent(ii) $d\mathbb{P}\times dt\text{-}a.e.,(\omega,t)\in \Omega\times [0,T]$, there exist constants $\gamma>0,\lambda>0,\gamma_{0}>0,M_{2}>0$ and an $\mathbb{F}$-progressively measurable process $\zeta=(\zeta_{t})$ with $\big\Vert
\int_{0}^{T}\zeta_{t}(\omega)dt\big\Vert_{\infty}\le M_{2}$, such that for all $z\in \mathbb{R}^{1\times d}$,
$$|f(t,z)|\le \zeta_{t}(\omega)+\psi(|U_{t}|)+\frac{\gamma}{2}|z|^{2}+n\lambda |V_{t}|^{1+\alpha}+\psi_{0}(\Vert U_{t}\Vert_{L^{2}(\Omega)})+\gamma_{0}\Vert V_{t}\Vert_{L^{2}(\Omega)}^{1+\alpha};$$

\noindent(iii) $d\mathbb{P}\times dt\text{-}a.e.,(\omega,t)\in \Omega\times [0,T]$, for all $z,\bar{z}\in \mathbb{R}^{1\times d}$,
$$|f(t,z)-f(t,\bar{z})|\le \psi(|U_{t}|\vee \Vert U_{t}\Vert_{L^{2}(\Omega)})(1+|z|+|\bar{z}|+2\Vert V_{t}\Vert_{L^{2}(\Omega)})|z-\bar{z}|.$$
\end{proposition}

\noindent Then the quadratic growth BSDE
\begin{equation}\label{eq 3.2}\tag{3.2}
Y_{t}=\xi+\int_{t}^{T}f(s,Z_{s})ds-\int_{t}^{T}Z_{s}dW_{s},\ \ t\in[0,T],
\end{equation}
has a unique adapted solution $(Y,Z)\in \mathcal{S}^{\infty}_{\mathbb{F}}(0,T;\mathbb{R})\times \mathcal{Z}^{2}_{\mathbb{F}}(0,T;\mathbb{R}^{1\times d})$. Moreover, for all $t\in[0,T]$ and stopping time $\tau\in \mathscr{T}[t,T]$, we have
\begin{equation*}\label{eq 3.3}\tag{3.3}
\begin{aligned}
|Y_{t}|&\le\ \frac{1}{\gamma}\ln 2+\Vert\xi\Vert_{\infty}+\bigg\Vert\int_{0}^{T}\zeta_{s}(\omega)ds\bigg\Vert_{\infty}+(\psi+\psi_{0})(\Vert U\Vert_{\mathcal{S}^{\infty}_{\mathbb{F}}(t,T)})(T-t)\\
&\hspace{2em}+\gamma_{0}\Vert V \Vert_{\mathcal{Z}^{2}_{\mathbb{F}}(t,T)}^{1+\alpha}(T-t)^{\frac{1-\alpha}{2}}+\gamma^{\frac{1+\alpha}{1-\alpha}}M_{n,\lambda,\alpha}\Vert V \Vert_{\mathcal{Z}^{2}_{\mathbb{F}}(t,T)}^{2\frac{1+\alpha}{1-\alpha}}(T-t)
\end{aligned}
\end{equation*}
and
\begin{equation*}\label{eq 3.4}\tag{3.4}
\begin{aligned}
\mathbb{E}_{\tau}\Big[\int_{\tau}^{T}|Z_{s}|^{2}ds\Big]\le \frac{1}{\gamma}\exp\Big\{2\gamma \Vert Y \Vert_{\mathcal{S}^{\infty}_{\mathbb{F}}(t,T)}\Big\}\bigg(2\bigg\Vert \int_{0}^{T}\zeta_{s}(\omega)ds\bigg\Vert_{\infty}+2(\psi+\psi_{0})(\Vert U\Vert_{\mathcal{S}^{\infty}_{\mathbb{F}}(t,T)})(T-t)\\
+1+2\gamma_{0}\Vert V \Vert_{\mathcal{Z}^{2}_{\mathbb{F}}(t,T)}^{1+\alpha}(T-t)^{\frac{1-\alpha}{2}}+2M_{n,\lambda,\alpha}\Vert V \Vert_{\mathcal{Z}^{2}_{\mathbb{F}}(t,T)}^{2\frac{1+\alpha}{1-\alpha}}(T-t)\bigg)+\frac{1}{\gamma^{2}}\exp\{2\gamma\Vert\xi\Vert_{\infty}\},
\end{aligned}
\end{equation*}
where
\begin{equation}\label{eq 3.5}\tag{3.5}
M_{n,\lambda,\alpha}:=\frac{1-\alpha}{2}(1+\alpha)^{\frac{1+\alpha}{1-\alpha}}(n\lambda)^{\frac{2}{1-\alpha}}.
\end{equation}

\begin{proof} On the one hand, using Young's inequality, we have for every $p\ge 1$ and $t\in [0,T]$,
\begin{equation*}\label{eq 3.6}\tag{3.6}
\begin{aligned}
p\gamma n\lambda |V_{s}|^{1+\alpha}
\le \frac{1}{2\Vert V \Vert_{\mathcal{Z}^{2}_{\mathbb{F}}(t,T)}^{2}}|V_{s}|^{2}+(p\gamma)^{\frac{2}{1-\alpha}} M_{n,\lambda,\alpha}\Vert V \Vert_{\mathcal{Z}^{2}_{\mathbb{F}}(t,T)}^{2\frac{1+\alpha}{1-\alpha}},\ \ 0\le t\le s\le T,
\end{aligned}
\end{equation*}
where $M_{n,\lambda,\alpha}$ is defined in (\ref{eq 3.5}). Then it can be inferred from the John-Nirenberg inequality for BMO martingales in Lemma 2.1 that, for all $\tau\in \mathscr{T}[t,T]$,
\begin{equation*}\label{eq 3.7}\tag{3.7}
\mathbb{E}_{\tau}\bigg[\exp\bigg(\frac{1}{2\Vert V \Vert_{\mathcal{Z}^{2}_{\mathbb{F}}(t,T)}^{2}}\int_{\tau}^{T}|V_{s}|^{2}ds\bigg)\bigg]\le 2.
\end{equation*}

\noindent On the other hand, applying H\"{o}lder's inequality yields that, for $t\in[0,T]$,
\begin{equation*}\label{eq 3.8}\tag{3.8}
\int_{t}^{T}(\mathbb{E}[|V_{s}|^{2}])^{\frac{1+\alpha}{2}}ds\le\Big(\int_{t}^{T}\mathbb{E}[|V_{s}|^{2}]ds\Big)^{\frac{1+\alpha}{2}}(T-t)^{\frac{1-\alpha}{2}}\le \Vert V \Vert_{\mathcal{Z}^{2}_{\mathbb{F}}(t,T)}^{1+\alpha}(T-t)^{\frac{1-\alpha}{2}}.
\end{equation*}

\noindent Combining this with (\ref{eq 3.6}), (\ref{eq 3.7}) and (\ref{eq 3.8}), we get
\begin{equation*}
\begin{aligned}
&\mathbb{E}_{t}\bigg[\exp\bigg\{p\gamma|\xi|+p\gamma\int_{t}^{T}\Big(\zeta_{s}(\omega)+\psi(|U_{s}|)+n\lambda|V_{s}|^{1+\alpha}+\psi_{0}(\Vert U_{s}\Vert_{L^{2}(\Omega)})+\gamma_{0}\Vert V_{s}\Vert^{1+\alpha}_{L^{2}(\Omega)}\Big)ds\bigg\}\bigg]\\
\le\ &2\exp\bigg\{p\gamma\Vert\xi\Vert_{\infty}+p\gamma\bigg\Vert\int_{0}^{T}\zeta_{s}(\omega)ds\bigg\Vert_{\infty}+p\gamma\psi(\Vert U\Vert_{\mathcal{S}^{\infty}_{\mathbb{F}}(t,T)})(T-t)+p\gamma\psi_{0}(\Vert U\Vert_{\mathcal{S}^{\infty}_{\mathbb{F}}(t,T)})(T-t)\\
&\hspace{3em}+p\gamma\gamma_{0}\Vert V \Vert_{\mathcal{Z}^{2}_{\mathbb{F}}(t,T)}^{1+\alpha}(T-t)^{\frac{1-\alpha}{2}}+(p\gamma)^{\frac{2}{1-\alpha}}M_{n,\lambda,\alpha}\Vert V \Vert_{\mathcal{Z}^{2}_{\mathbb{F}}(t,T)}^{2\frac{1+\alpha}{1-\alpha}}(T-t)\bigg\}<+\infty.
\end{aligned}
\end{equation*}

\noindent Thanks to the conditions (i)-(iii), we can use Proposition 3 and Corollary 4 in Briand and Hu \cite{BH08} to conclude that BSDE (\ref{eq 3.2}) admits a solution $(Y,Z)$ such that $\mathbb{E}\big[\int_{0}^{T}|Z_{s}|^{2}ds\big]<+\infty$ and
\begin{equation*}
\begin{aligned}
\exp\{\gamma |Y_{t}|\}&\le\ 2\exp\bigg\{\gamma\Vert\xi\Vert_{\infty}+\gamma\bigg\Vert\int_{0}^{T}\zeta_{s}(\omega)ds\bigg\Vert_{\infty}+\gamma(\psi+\psi_{0})(\Vert U\Vert_{\mathcal{S}^{\infty}_{\mathbb{F}}(t,T)})(T-t)\\
&\hspace{2em}+\gamma\gamma_{0}\Vert V \Vert_{\mathcal{Z}^{2}_{\mathbb{F}}(t,T)}^{1+\alpha}(T-t)^{\frac{1-\alpha}{2}}+\gamma^{\frac{2}{1-\alpha}}M_{n,\lambda,\alpha}\Vert V \Vert_{\mathcal{Z}^{2}_{\mathbb{F}}(t,T)}^{2\frac{1+\alpha}{1-\alpha}}(T-t)\bigg\}.
\end{aligned}
\end{equation*}
This implies that (\ref{eq 3.3}) holds and $Y\in \mathcal{S}^{\infty}_{\mathbb{F}}(0,T;\mathbb{R})$.\vspace{-0.5em}\\

To prove that $Z\cdot W$ is a BMO martingale, we apply It\^{o}-Tanaka's formula to $\exp(2\gamma |Y_{t}|)$. This yields for all $t\in [0,T]$ and stopping time $\tau\in \mathscr{T}[t,T]$,
\begin{equation*}\label{eq 3.9}\tag{3.9}
\begin{aligned}
\exp\{2\gamma |Y_{\tau}|\}+2\gamma^{2}\mathbb{E}_{\tau}\Big[\int_{\tau}^{T}\exp\{2\gamma|Y_{s}|\}|Z_{s}|^{2}ds\Big]\le \mathbb{E}_{\tau}\big[\exp\{2\gamma|\xi|\}\big]+2\gamma\mathbb{E}_{\tau}\Big[\int_{\tau}^{T}\exp\{2\gamma |Y_{s}|\}\\
\cdot\Big(\zeta_{s}(\omega)+\psi(|U_{s}|)+n\lambda|V_{s}|^{1+\alpha}+\frac{\gamma}{2}|Z_{s}|^{2}+\psi_{0}(\Vert U_{s} \Vert_{L^{2}(\Omega)})+\gamma_{0}\Vert V_{s} \Vert_{L^{2}(\Omega)}^{1+\alpha}\Big)ds\Big].
\end{aligned}
\end{equation*}

\noindent Consequently, we have that
\begin{equation*}
\begin{aligned}
\gamma^{2}\mathbb{E}_{\tau}\Big[\int_{\tau}^{T}|Z_{s}|^{2}ds\Big]\le 2\gamma\exp\Big\{2\gamma \Vert Y \Vert_{\mathcal{S}^{\infty}_{\mathbb{F}}(t,T)}\Big\}\bigg(\bigg\Vert \int_{0}^{T}\zeta_{s}(\omega)ds\bigg\Vert_{\infty}+(\psi+\psi_{0})(\Vert U\Vert_{\mathcal{S}^{\infty}_{\mathbb{F}}(t,T)})(T-t)\\
+\frac{1}{2}+\gamma_{0}\Vert V \Vert_{\mathcal{Z}^{2}_{\mathbb{F}}(t,T)}^{1+\alpha}(T-t)^{\frac{1-\alpha}{2}}+M_{n,\lambda,\alpha}\Vert V \Vert_{\mathcal{Z}^{2}_{\mathbb{F}}(t,T)}^{2\frac{1+\alpha}{1-\alpha}}(T-t)\bigg)+\exp\{2\gamma\Vert\xi\Vert_{\infty}\}.
\end{aligned}
\end{equation*}
Relation (\ref{eq 3.4}) now follows easily.

Finally, by following the same argument as in the proof of Hu and Tang \cite[Lemma 2.1]{HT16}, we use the Girsanov transformation to obtain a comparison theorem for the solutions of BSDE (3.2), which implies the uniqueness. Thus, the proof is complete.
\end{proof}

In order to prove the existence and the uniqueness of mean-field BSDE (\ref{eq 3.1}) in a subset of $\mathcal{S}^{\infty}_{\mathbb{F}}(T-\varepsilon,T;\mathbb{R}^{n})\times \mathcal{Z}^{2}_{\mathbb{F}}(T-\varepsilon,T;\mathbb{R}^{n\times d})$, it is necessary to define the following set, with positive constants $K_{1},K_{2},\varepsilon$, assigned accordingly:
\begin{equation*}\label{eq 3.10}\tag{3.10}
\begin{aligned}
&\mathcal{B}_{\varepsilon}(K_{1},K_{2}):=\Big\{(Y,Z)\in \mathcal{S}^{\infty}_{\mathbb{F}}(T-\varepsilon,T;\mathbb{R}^{n})\times \mathcal{Z}^{2}_{\mathbb{F}}(T-\varepsilon,T;\mathbb{R}^{n\times d}):\\
&\hspace{15.6em}\Vert Y\Vert_{\mathcal{S}^{\infty}_{\mathbb{F}}(T-\varepsilon,T)}\le K_{1}\ \ \text{and}\ \ \Vert Z\Vert^{2}_{\mathcal{Z}^{2}_{\mathbb{F}}(T-\varepsilon,T)}\le K_{2}\Big\}
\end{aligned}
\end{equation*}
equipped with the norm
$$\Vert (Y,Z)\Vert_{\mathcal{B}_{\varepsilon}(K_{1},K_{2})}:=\Big\{ \Vert Y\Vert^{2}_{\mathcal{S}^{\infty}_{\mathbb{F}}(T-\varepsilon,T)}+\Vert Z\Vert^{2}_{\mathcal{Z}^{2}_{\mathbb{F}}(T-\varepsilon,T)}\Big\}^{\frac{1}{2}}.$$

\begin{theorem}
Let the assumptions (A1)-(A3) be satisfied. Then there exists $\varepsilon>0$, depending only on the constants $n, \gamma,\gamma_{0},\lambda,\alpha,M_{1},M_{2}$ and the functions $\psi,\ \psi_{0}$, such that on the time interval $[T-\varepsilon,T]$, the mean-field BSDE (\ref{eq 3.1}) admits a unique local solution $(Y,Z)\in \mathcal{B}_{\varepsilon}(K_{1},K_{2})$ with
\begin{equation}\label{eq 3.11}\tag{3.11}
K_{1}=\frac{2n}{\gamma}\ln2+2n(M_{1}+M_{2}),\ \ K_{2}=\frac{2n}{\gamma^{2}}e^{2\gamma M_{1}}+\frac{2n}{\gamma}e^{2\gamma K_{1}}(1+2M_{2}).
\end{equation}
\end{theorem}

\begin{proof} For $1\le i\le n,\ H\in \mathbb{R}^{n\times d}$ and $z\in \mathbb{R}^{1\times d}$, we denote by $H(z;i)$ the matrix in $\mathbb{R}^{n\times d}$ whose $i$th row is $z$ and whose $j$th row is $H^{j}$ for any $j\ne i$.

\noindent\textbf{Step 1.}\ \ Letting $(U,V)\in \mathcal{S}^{\infty}_{\mathbb{F}}(0,T;\mathbb{R}^{n})\times \mathcal{Z}^{2}_{\mathbb{F}}(0,T;\mathbb{R}^{n\times d})$, we consider the following decoupled system of quadratic BSDEs:
\begin{equation}\label{eq 3.12}\tag{3.12}
Y_{t}^{i}=\xi^{i}+\int_{t}^{T}f^{i}(s,U_{s},V_{s}(Z_{s}^{i};i),\mathbb{P}_{(U_{s},V_{s})})ds-\int_{t}^{T}Z_{s}^{i}dW_{s},\ \ t\in[0,T],\ \ 1\le i\le n.
\end{equation}

\noindent For every $1\le i\le n$, we denote by $f^{i,U,V}$ the generator:
$$f^{i,U,V}(t,z)=f^{i}(t,U_{t},V_{t}(z;i),\mathbb{P}_{(U_{t},V_{t})}),\ \ (t,z)\in [0,T]\times \mathbb{R}^{1\times d}.$$

\noindent Based on the assumptions (A2)-(A3), we verify that, $d\mathbb{P}\times dt\text{-}a.e.$, for $z,\bar{z}\in \mathbb{R}^{1\times d}$,
$$|f^{i,U,V}(t,z)|\le \zeta_{t}(\omega)+\psi(|U_{t}|)+\frac{\gamma}{2}|z|^{2}+n\lambda |V_{t}|^{1+\alpha}+\psi_{0}(\Vert U_{t}\Vert_{L^{2}(\Omega)})+\gamma_{0}\Vert V_{t}\Vert_{L^{2}(\Omega)}$$
and
$$|f^{i,U,V}(t,z)-f^{i,U,V}(t,\bar{z})|\le \psi(|U_{t}|\vee \Vert U_{t}\Vert_{L^{2}(\Omega)})(1+|z|+|\bar{z}|+2\Vert V_{t}\Vert_{L^{2}(\Omega)})|z-\bar{z}|.$$

\noindent Then in view of assumption (A1), it follows from Proposition 3.1 that BSDE (\ref{eq 3.12}) possesses a unique solution $(Y,Z)\in \mathcal{S}_{\mathbb{F}}^{\infty}(0,T;\mathbb{R}^{n})\times \mathcal{Z}_{\mathbb{F}}^{2}(0,T;\mathbb{R}^{n\times d})$. Therefore, we can define a mapping $\Psi$ on $\mathcal{S}^{\infty}_{\mathbb{F}}(0,T;\mathbb{R}^{n})\times \mathcal{Z}^{2}_{\mathbb{F}}(0,T;\mathbb{R}^{n\times d})$ as follows:
\begin{equation*}
\Psi(U,V):=(Y,Z),\ \ (U,V)\in \mathcal{S}^{\infty}_{\mathbb{F}}(0,T;\mathbb{R}^{n})\times \mathcal{Z}^{2}_{\mathbb{F}}(0,T;\mathbb{R}^{n\times d}).
\end{equation*}

\noindent\textbf{Step 2.}\ \ We prove that, for suitably chosen $K_{1},K_{2}>0$, the mapping $\Psi$ maps $\mathcal{B}_{\varepsilon}(K_{1},K_{2})$ into itself, where $\mathcal{B}_{\varepsilon}(K_{1},K_{2})$ is defined in (\ref{eq 3.10}). We observe that by Proposition 3.1, for every $1\le i\le n,\ t\in[0,T]$ and stopping time $\tau\in \mathscr{T}[t,T]$, we have the following estimates:
\begin{equation*}\label{eq 3.13}\tag{3.13}
\begin{aligned}
|Y_{t}^{i}|&\le\ \frac{1}{\gamma}\ln 2+\Vert\xi^{i}\Vert_{\infty}+\bigg\Vert\int_{0}^{T}\zeta_{s}(\omega)ds\bigg\Vert_{\infty}+(\psi+\psi_{0})(\Vert U\Vert_{\mathcal{S}^{\infty}_{\mathbb{F}}(t,T)})(T-t)\\
&\hspace{2em}+\gamma_{0}\Vert V \Vert_{\mathcal{Z}^{2}_{\mathbb{F}}(t,T)}^{1+\alpha}(T-t)^{\frac{1-\alpha}{2}}+\gamma^{\frac{1+\alpha}{1-\alpha}}M_{n,\lambda,\alpha}\Vert V \Vert_{\mathcal{Z}^{2}_{\mathbb{F}}(t,T)}^{2\frac{1+\alpha}{1-\alpha}}(T-t)
\end{aligned}
\end{equation*}
and
\begin{equation*}\label{eq 3.14}\tag{3.14}
\begin{aligned}
\mathbb{E}_{\tau}\Big[\int_{\tau}^{T}|Z_{s}^{i}|^{2}ds\Big]\le \frac{1}{\gamma}\exp\Big\{2\gamma \Vert Y^{i} \Vert_{\mathcal{S}^{\infty}_{\mathbb{F}}(t,T)}\Big\}\bigg(2\bigg\Vert \int_{0}^{T}\zeta_{s}(\omega)ds\bigg\Vert_{\infty}+2(\psi+\psi_{0})(\Vert U\Vert_{\mathcal{S}^{\infty}_{\mathbb{F}}(t,T)})(T-t)\\
+1+2\gamma_{0}\Vert V \Vert_{\mathcal{Z}^{2}_{\mathbb{F}}(t,T)}^{1+\alpha}(T-t)^{\frac{1-\alpha}{2}}+2M_{n,\lambda,\alpha}\Vert V \Vert_{\mathcal{Z}^{2}_{\mathbb{F}}(t,T)}^{2\frac{1+\alpha}{1-\alpha}}(T-t)\bigg)+\frac{1}{\gamma^{2}}\exp\{2\gamma\Vert\xi^{i}\Vert_{\infty}\},
\end{aligned}
\end{equation*}
where $M_{n,\lambda,\alpha}$ is defined in (\ref{eq 3.5}). Then, in view of (A1)-(A2), we get that, for all $t\in[0,T]$,
\begin{equation*}\label{eq 3.15}\tag{3.15}
\begin{aligned}
\Vert Y \Vert_{S^{\infty}_{\mathbb{F}}(t,T)}&\le\ \frac{n}{\gamma}\ln 2+n(M_{1}+M_{2})+n(\psi+\psi_{0})(\Vert U\Vert_{\mathcal{S}^{\infty}_{\mathbb{F}}(t,T)})(T-t)\\
&\hspace{2em}+n\gamma_{0}\Vert V \Vert_{\mathcal{Z}^{2}_{\mathbb{F}}(t,T)}^{1+\alpha}(T-t)^{\frac{1-\alpha}{2}}+n\gamma^{\frac{1+\alpha}{1-\alpha}}M_{n,\lambda,\alpha}\Vert V \Vert_{\mathcal{Z}^{2}_{\mathbb{F}}(t,T)}^{2\frac{1+\alpha}{1-\alpha}}(T-t)
\end{aligned}
\end{equation*}
and
\begin{equation*}\label{eq 3.16}\tag{3.16}
\begin{aligned}
\Vert Z\Vert_{\mathcal{Z}^{2}_{\mathbb{F}}(t,T)}^{2}\le \frac{n}{\gamma}\exp\Big\{2\gamma \Vert Y \Vert_{\mathcal{S}^{\infty}_{\mathbb{F}}(t,T)}\Big\}\bigg(1+2M_{2}+2(\psi+\psi_{0})(\Vert U\Vert_{\mathcal{S}^{\infty}_{\mathbb{F}}(t,T)})(T-t)\\
+2\gamma_{0}\Vert V \Vert_{\mathcal{Z}^{2}_{\mathbb{F}}(t,T)}^{1+\alpha}(T-t)^{\frac{1-\alpha}{2}}+2M_{n,\lambda,\alpha}\Vert V \Vert_{\mathcal{Z}^{2}_{\mathbb{F}}(t,T)}^{2\frac{1+\alpha}{1-\alpha}}(T-t)\bigg)+\frac{n}{\gamma^{2}}\exp\{2\gamma M_{1}\}.
\end{aligned}
\end{equation*}
Let us define
\begin{equation}\label{eq 3.17}\tag{3.17}
K_{1}:=\frac{2n}{\gamma}\ln2+2n(M_{1}+M_{2}),\ \ K_{2}:=\frac{2n}{\gamma^{2}}e^{2\gamma M_{1}}+\frac{2n}{\gamma}e^{2\gamma K_{1}}(1+2M_{2}).
\end{equation}
We denote by $x_{1}$ and $x_{2}$ the unique solutions to the following equations:
\begin{equation*}\label{eq 3.18}\tag{3.18}
n(\psi+\psi_{0})(K_{1})x_{1}+n\gamma_{0}K_{2}^{\frac{1+\alpha}{2}}x_{1}^{\frac{1-\alpha}{2}}+n\gamma^{\frac{1+\alpha}{1-\alpha}}M_{n,\lambda,\alpha}K_{2}^{\frac{1+\alpha}{1-\alpha}}x_{1}=\frac{K_{1}}{2}
\end{equation*}
and
\begin{equation*}\label{eq 3.19}\tag{3.19}
2(\psi+\psi_{0})(K_{1})x_{2}+2\gamma_{0}K_{2}^{\frac{1+\alpha}{2}}x_{2}^{\frac{1-\alpha}{2}}+2M_{n,\lambda,\alpha}K_{2}^{\frac{1+\alpha}{1-\alpha}}x_{2}=\frac{\gamma K_{2}}{2n}e^{-2\gamma K_{1}},
\end{equation*}
respectively. By virtue of (\ref{eq 3.15})-(\ref{eq 3.19}), we check that
\begin{equation*}
\Vert U\Vert_{\mathcal{S}^{\infty}_{\mathbb{F}}(T-\varepsilon,T)}\le K_{1}\ \ \text{and}\ \ \Vert V \Vert_{\mathcal{Z}^{2}_{\mathbb{F}}(T-\varepsilon,T)}^{2}\le K_{2},
\end{equation*}
implies
\begin{equation*}
\Vert Y\Vert_{\mathcal{S}^{\infty}_{\mathbb{F}}(T-\varepsilon,T)}\le K_{1}\ \ \text{and}\ \ \Vert Z \Vert_{\mathcal{Z}^{2}_{\mathbb{F}}(T-\varepsilon,T)}^{2}\le K_{2},\ \ \text{for all}\ \varepsilon\in(0,x_{0}],
\end{equation*}
where $x_{0}=x_{1}\wedge x_{2}>0$. This means that
\begin{equation*}
\Psi(U,V)\in \mathcal{B}_{\varepsilon}(K_{1},K_{2}),\ \ \text{for all}\ (U,V)\in \mathcal{B}_{\varepsilon}(K_{1},K_{2}).
\end{equation*}
That is, the mapping $\Psi$ is stable in $\mathcal{B}_{\varepsilon}(K_{1},K_{2})$.

\noindent\textbf{Step 3.}\ \ We show that the mapping $\Psi$ is contractive in $\mathcal{B}_{\varepsilon}(K_{1},K_{2})$.

Indeed, for any fixed $\varepsilon\in(0,x_{0}]$, and $(U,V),(\widetilde{U},\widetilde{V})\in \mathcal{B}_{\varepsilon}(K_{1},K_{2})$, we set
\begin{equation*}
(Y,Z):=\Psi(U,V),\ \ (\widetilde{Y},\widetilde{Z}):=\Psi(\widetilde{U},\widetilde{V}).
\end{equation*}
Moreover, for all $t\in [T-\varepsilon,T]$, we put
\begin{equation*}
\Delta Y_{t}:=Y_{t}-\widetilde{Y}_{t},\ \ \Delta Z_{t}:=Z_{t}-\widetilde{Z}_{t},\ \ \Delta U_{t}:=U_{t}-\widetilde{U}_{t},\ \ \Delta V_{t}:=V_{t}-\widetilde{V}_{t}.
\end{equation*}
Then, for $1\le i\le n$ and $t\in[T-\varepsilon,T]$, we have
\begin{equation*}\label{eq 3.20}\tag{3.20}
\Delta Y_{t}^{i}=\int_{t}^{T}(I^{1,i}_{s}+I^{2,i}_{s})ds-\int_{t}^{T}\Delta Z_{s}^{i}dW_{s},
\end{equation*}
where
\begin{equation*}
I^{1,i}_{s}:=f^{i}(s,U_{s},V_{s}(Z_{s}^{i};i),\mathbb{P}_{(U_{s},V_{s})})-f^{i}(s,U_{s},V_{s}(\widetilde{Z}_{s}^{i};i),\mathbb{P}_{(U_{s},V_{s})})
\end{equation*}
and
\begin{equation*}
I^{2,i}_{s}:=f^{i}(s,U_{s},V_{s}(\widetilde{Z}_{s}^{i};i),\mathbb{P}_{(U_{s},V_{s})})-f^{i}(s,\widetilde{U}_{s},\widetilde{V}_{s}(\widetilde{Z}_{s}^{i};i),\mathbb{P}_{(\widetilde{U}_{s},\widetilde{V}_{s})}).
\end{equation*}
We notice that for $1\le i\le n$, $I_{s}^{1,i}$ can be written as $I_{s}^{1,i}=(Z_{s}^{i}-\widetilde{Z}_{s}^{i})\Lambda_{s}^{i}$, where
\begin{equation*}
\Lambda_{s}^{i}=\left\{
\begin{aligned}
&\frac{\big(f^{i}(s,U_{s},V_{s}(Z_{s}^{i};i),\mathbb{P}_{(U_{s},V_{s})})-f^{i}(s,U_{s},V_{s}(\widetilde{Z}_{s}^{i};i),\mathbb{P}_{(U_{s},V_{s})})\big)(Z_{s}^{i}-\widetilde{Z}_{s}^{i})^{\mathsf{T}}}{|Z_{s}^{i}-\widetilde{Z}_{s}^{i}|^{2}},\ \ \ \text{if}\ Z_{s}^{i}\ne\widetilde{Z}_{s}^{i};\\
&0,\hspace{30.45em}\text{if}\ Z_{s}^{i}=\widetilde{Z}_{s}^{i}.
\end{aligned}
\right.
\end{equation*}
By assumption (A3) we have, for every $s\in[t,T]$,
\begin{equation*}
|\Lambda_{s}^{i}|\le \psi(|U_{s}|+\Vert U_{s}\Vert_{L^{2}(\Omega)})\big(1+2|V_{s}|+|Z_{s}|+|\widetilde{Z}_{s}|+2\Vert V_{s}\Vert_{L^{2}(\Omega)}\big).
\end{equation*}
For $1\le i\le n$, we define
\begin{equation*}\label{eq 3.21}\tag{3.21}
\widetilde{W}_{t}^{i}:=W_{t}-\int_{0}^{t}\Lambda_{s}^{i}ds\ \ \text{and}\ \ d\mathbb{Q}^{i}:=\mathscr{E}(\Lambda^{i}\cdot W)_{0}^{T}d\mathbb{P}.
\end{equation*}
Notice that $\mathbb{E}^{\mathbb{P}}[\mathscr{E}(\Lambda^{i}\cdot W)_{0}^{T}]=1$ (refer to \cite[Theorem 2.3]{K94}). Then $\mathbb{Q}^{i}$ is equivalent to $\mathbb{P}$ and the process $\widetilde{W}^{i}$ is a Brownian motion under $\mathbb{Q}^{i}$. Therefore, we can reformulate BSDE (\ref{eq 3.20}) as:
\begin{equation*}\label{eq 3.22}\tag{3.22}
\Delta Y_{t}^{i}+\int_{t}^{T}\Delta Z_{s}^{i}d\widetilde{W}^{i}_{s}=\int_{t}^{T}I^{2,i}_{s}ds.
\end{equation*}
Taking the square and then the conditional expectation with respect to $\mathbb{Q}^{i}$ on both sides of the latter equation, we deduce that, for $1\le i\le n$,
\begin{equation*}\label{eq 3.23}\tag{3.23}
|\Delta Y_{t}^{i}|^{2}+\mathbb{E}_{t}^{\mathbb{Q}^{i}}\Big[\int_{t}^{T}|\Delta Z_{s}^{i}|^{2}ds\Big]=\mathbb{E}_{t}^{\mathbb{Q}^{i}}\bigg[\Big(\int_{t}^{T}I^{2,i}_{s}ds\Big)^{2}\bigg].
\end{equation*}
Following assumption (A3), we see that
\begin{equation*}
\begin{aligned}
&|I_{s}^{2,i}|\le
\bigg[\psi(K_{1})\Big(1+|V_{s}|+|\widetilde{V}_{s}|+2|\widetilde{Z}_{s}|+\Vert V_{s}\Vert_{L^{2}(\Omega)}+\Vert \widetilde{V}_{s}\Vert_{L^{2}(\Omega)}\Big)2\Vert \Delta U\Vert_{\mathcal{S}^{\infty}_{\mathbb{F}}(t,T)}\\
&\hspace{3.1em}+\Big(1+|V_{s}|^{\alpha}+|\widetilde{V}_{s}|^{\alpha}+2|\widetilde{Z}_{s}|^{\alpha}+\Vert V_{s}\Vert_{L^{2}(\Omega)}^{\alpha}+\Vert \widetilde{V}_{s}\Vert_{L^{2}(\Omega)}^{\alpha}\Big)\Big(\sqrt{n}|\Delta V_{s}|+\Vert \Delta V_{s}\Vert_{L^{2}(\Omega)}\Big)\bigg].\\
\end{aligned}
\end{equation*}

\noindent Thus, from (\ref{eq 3.23}) we derive that
\begin{equation*}\label{eq 3.24}\tag{3.24}
\begin{aligned}
&\hspace{2.6em}|\Delta Y_{t}^{i}|^{2}+\mathbb{E}_{t}^{\mathbb{Q}^{i}}\Big[\int_{t}^{T}|\Delta Z_{s}^{i}|^{2}ds\Big]\\ &\le12\psi(K_{1})^{2}\mathbb{E}_{t}^{\mathbb{Q}^{i}}\Big[\Big(\int_{t}^{T}\big(1+|V_{s}|+|\widetilde{V}_{s}|+2|\widetilde{Z}_{s}|+\Vert V_{s}\Vert_{L^{2}(\Omega)}+\Vert \widetilde{V}_{s}\Vert_{L^{2}(\Omega)}\big)\Vert \Delta U\Vert_{\mathcal{S}^{\infty}_{\mathbb{F}}(t,T)}ds\Big)^{2}\Big]\\
&\hspace{1.5em}+3n\psi(K_{1})^{2}\mathbb{E}_{t}^{\mathbb{Q}^{i}}\Big[\Big(\int_{t}^{T}\big(1+|V_{s}|^{\alpha}+|\widetilde{V}_{s}|^{\alpha}+2|\widetilde{Z}_{s}|^{\alpha}+\Vert V_{s}\Vert_{L^{2}(\Omega)}^{\alpha}+\Vert \widetilde{V}_{s}\Vert_{L^{2}(\Omega)}^{\alpha}\big)|\Delta V_{s}|ds\Big)^{2}\Big]\\
&\hspace{1.5em}+3\psi(K_{1})^{2}\mathbb{E}_{t}^{\mathbb{Q}^{i}}\Big[\Big(\int_{t}^{T}\big(1+|V_{s}|^{\alpha}+|\widetilde{V}_{s}|^{\alpha}+2|\widetilde{Z}_{s}|^{\alpha}+\Vert V_{s}\Vert_{L^{2}(\Omega)}^{\alpha}+\Vert \widetilde{V}_{s}\Vert_{L^{2}(\Omega)}^{\alpha}\big)\Vert \Delta V_{s}\Vert_{L^{2}(\Omega)}ds\Big)^{2}\Big].
\end{aligned}
\end{equation*}
We now handle the first term in the right hand side of the last inequality, using H\"{o}lder's inequality and Lemma 2.3 implies that, for $t\in[T-\varepsilon,T]$,
\begin{equation*}\label{eq 3.25}\tag{3.25}
\begin{aligned}
&\hspace{1.5em}12\psi(K_{1})^{2}\mathbb{E}_{t}^{\mathbb{Q}^{i}}\Big[\Big(\int_{t}^{T}\big(1+|V_{s}|+|\widetilde{V}_{s}|+2|\widetilde{Z}_{s}|+\Vert V_{s}\Vert_{L^{2}(\Omega)}+\Vert \widetilde{V}_{s}\Vert_{L^{2}(\Omega)}\big)\Vert \Delta U\Vert_{\mathcal{S}^{\infty}_{\mathbb{F}}(t,T)}ds\Big)^{2}\Big]\\
&\le72\varepsilon\psi(K_{1})^{2}\Vert\Delta U\Vert_{\mathcal{S}^{\infty}_{\mathbb{F}}(t,T)}^{2}\bigg(T+\mathbb{E}_{t}^{\mathbb{Q}^{i}}\Big[\int_{t}^{T}\big(|V_{s}|^{2}+|\widetilde{V}_{s}|^{2}+4|\widetilde{Z}_{s}|^{2}\big)ds\Big]+\int_{t}^{T}\mathbb{E}^{\mathbb{P}}\Big[|V_{s}|^{2}+|\widetilde{V}_{s}|^{2}\Big]ds\bigg)\\
&\le72\varepsilon\psi(K_{1})^{2}(T+6c_{2}^{2}K_{2}+2K_{2})\Vert\Delta U\Vert_{\mathcal{S}^{\infty}_{\mathbb{F}}(t,T)}^{2}.
\end{aligned}
\end{equation*}
For the second term, we deduce from H\"{o}lder's inequality that
\begin{equation*}
\begin{aligned}
&\hspace{1.5em}\mathbb{E}_{t}^{\mathbb{Q}^{i}}\Big[\Big(\int_{t}^{T}\big(1+|V_{s}|^{\alpha}+|\widetilde{V}_{s}|^{\alpha}+2|\widetilde{Z}_{s}|^{\alpha}+\Vert V_{s}\Vert_{L^{2}(\Omega)}^{\alpha}+\Vert \widetilde{V}_{s}\Vert_{L^{2}(\Omega)}^{\alpha}\big)|\Delta V_{s}|ds\Big)^{2}\Big]\\
&\le\mathbb{E}_{t}^{\mathbb{Q}^{i}}\Big[\int_{t}^{T}\Big(1+|V_{s}|^{\alpha}+|\widetilde{V}_{s}|^{\alpha}+2|\widetilde{Z}_{s}|^{\alpha}+\Vert V_{s}\Vert_{L^{2}(\Omega)}^{\alpha}+\Vert \widetilde{V}_{s}\Vert_{L^{2}(\Omega)}^{\alpha}\Big)^{2}ds\cdot\int_{t}^{T}|\Delta V_{s}|^{2}ds\Big]\\
&\le6\mathbb{E}_{t}^{\mathbb{Q}^{i}}\Big[\Big(\int_{t}^{T}\big(1+|V_{s}|^{2\alpha}+|\widetilde{V}_{s}|^{2\alpha}+4|\widetilde{Z}_{s}|^{2\alpha}+\Vert V_{s}\Vert_{L^{2}(\Omega)}^{2\alpha}+\Vert \widetilde{V}_{s}\Vert_{L^{2}(\Omega)}^{2\alpha}\big)ds\Big)^{2}\Big]^{\frac{1}{2}}\mathbb{E}_{t}^{\mathbb{Q}^{i}}\Big[\Big(\int_{t}^{T}|\Delta V_{s}|^{2}ds\Big)^{2}\Big]^{\frac{1}{2}}.
\end{aligned}
\end{equation*}

\noindent By using H\"{o}lder's inequality and Young's inequality, as well as the Lemmas 2.2 and 2.3, we get, for $t\in[T-\varepsilon,T]$, that
\begin{equation*}
\begin{aligned}
&\hspace{1.5em}\mathbb{E}_{t}^{\mathbb{Q}^{i}}\Big[\Big(\int_{t}^{T}\big(1+|V_{s}|^{2\alpha}+|\widetilde{V}_{s}|^{2\alpha}+4|\widetilde{Z}_{s}|^{2\alpha}+\Vert V_{s}\Vert_{L^{2}(\Omega)}^{2\alpha}+\Vert \widetilde{V}_{s}\Vert_{L^{2}(\Omega)}^{2\alpha}\big)ds\Big)^{2}\Big]^{\frac{1}{2}}\\
&\le \varepsilon^{1-\alpha}\mathbb{E}_{t}^{\mathbb{Q}^{i}}\Big[\Big(T^{\alpha}+8(1-\alpha)+\alpha\int_{t}^{T}\big(|V_{s}|^{2}+|\widetilde{V}_{s}|^{2}+4|\widetilde{Z}_{s}|^{2}\big)ds+\alpha\int_{t}^{T}\mathbb{E}^{\mathbb{P}}[|V_{s}|^{2}+|\widetilde{V}_{s}|^{2}]ds\Big)^{2}\Big]^{\frac{1}{2}}\\
&\le\varepsilon^{1-\alpha}\Big[T^{\alpha}+8+\alpha L_{4}^{2}c_{2}^{2}\Vert V\cdot W\Vert_{BMO(\mathbb{P})}^{2}+\alpha L_{4}^{2}c_{2}^{2}\Vert \widetilde{V}\cdot W\Vert_{BMO(\mathbb{P})}^{2}+4\alpha L_{4}^{2}c_{2}^{2}\Vert \widetilde{Z}\cdot W\Vert_{BMO(\mathbb{P})}^{2}+2\alpha K_{2}\Big]\\
&\le\varepsilon^{1-\alpha}\Big(T^{\alpha}+8+6\alpha L_{4}^{2}c_{2}^{2}K_{2}+2\alpha K_{2}\Big).
\end{aligned}
\end{equation*}
Similarly, it follows from Lemma 2.2 and Lemma 2.3 that
\begin{equation*}
\begin{aligned}
\mathbb{E}_{t}^{\mathbb{Q}^{i}}\Big[\Big(\int_{t}^{T}|\Delta V_{s}|^{2}ds\Big)^{2}\Big]^{\frac{1}{2}}\le L_{4}^{2}\Vert \Delta V\cdot \widetilde{W}^{i}\Vert_{BMO(\mathbb{Q}^{i})}^{2}\le L_{4}^{2}c_{2}^{2}\Vert \Delta V\cdot W\Vert_{BMO(\mathbb{P})}^{2}= L_{4}^{2}c_{2}^{2}\Vert\Delta V\Vert_{\mathcal{Z}^{2}_{\mathbb{F}}(t,T)}^{2}.
\end{aligned}
\end{equation*}
Now combining the above estimates we obtain
\begin{equation*}\label{eq 3.26}\tag{3.26}
\begin{aligned}
&\hspace{1.5em}\mathbb{E}_{t}^{\mathbb{Q}^{i}}\Big[\Big(\int_{t}^{T}\big(1+|V_{s}|^{\alpha}+|\widetilde{V}_{s}|^{\alpha}+2|\widetilde{Z}_{s}|^{\alpha}+\Vert V_{s}\Vert_{L^{2}(\Omega)}^{\alpha}+\Vert \widetilde{V}_{s}\Vert_{L^{2}(\Omega)}^{\alpha}\big)|\Delta V_{s}|ds\Big)^{2}\Big]\\
&\le 6\varepsilon^{1-\alpha}L_{4}^{2}c_{2}^{2}\Big(T^{\alpha}+8+6\alpha L_{4}^{2}c_{2}^{2}K_{2}+2\alpha K_{2}\Big)\Vert\Delta V\Vert_{\mathcal{Z}^{2}_{\mathbb{F}}(t,T)}^{2}.
\end{aligned}
\end{equation*}
Note that $\alpha\in[0,1)$. Thus, using H\"{o}lder's inequality, and Fubini's theorem together with Lemma 2.3, we have
\begin{equation*}\label{eq 3.27}\tag{3.27}
\begin{aligned}
&\hspace{1.5em}3\psi(K_{1})^{2}\mathbb{E}_{t}^{\mathbb{Q}^{i}}\Big[\Big(\int_{t}^{T}\big(1+|V_{s}|^{\alpha}+|\widetilde{V}_{s}|^{\alpha}+2|\widetilde{Z}_{s}|^{\alpha}+\Vert V_{s}\Vert_{L^{2}(\Omega)}^{\alpha}+\Vert \widetilde{V}_{s}\Vert_{L^{2}(\Omega)}^{\alpha}\big)\Vert \Delta V_{s}\Vert_{L^{2}(\Omega)}ds\Big)^{2}\Big]\\
&\le18\psi(K_{1})^{2}\mathbb{E}_{t}^{\mathbb{Q}^{i}}\Big[\int_{t}^{T}\big(1+|V_{s}|^{2\alpha}+|\widetilde{V}_{s}|^{2\alpha}+4|\widetilde{Z}_{s}|^{2\alpha}+\Vert V_{s}\Vert_{L^{2}(\Omega)}^{2\alpha}+\Vert \widetilde{V}_{s}\Vert_{L^{2}(\Omega)}^{2\alpha}\big)ds\Big]\int_{t}^{T}\Vert \Delta V_{s}\Vert_{L^{2}(\Omega)}^{2}ds\\
&\le 18\varepsilon^{1-\alpha}\psi(K_{1})^{2}\bigg[T^{\alpha}+6-6\alpha+\alpha\mathbb{E}_{t}^{\mathbb{Q}^{i}}\Big[\int_{t}^{T}|V_{s}|^{2}ds\Big]+\alpha\mathbb{E}_{t}^{\mathbb{Q}^{i}}\Big[\int_{t}^{T}|\widetilde{V}_{s}|^{2}ds\Big]+4\alpha\mathbb{E}_{t}^{\mathbb{Q}^{i}}\Big[\int_{t}^{T}|\widetilde{Z}_{s}|^{2}ds\Big]\\
&\hspace{6em}+\Big(\mathbb{E}^{\mathbb{P}}\Big[\int_{t}^{T}|V_{s}|^{2}ds\Big]\Big)^{\alpha}+\Big(\mathbb{E}^{\mathbb{P}}\Big[\int_{t}^{T}|\widetilde{V}_{s}|^{2}ds\Big]\Big)^{\alpha}\bigg]\cdot\Vert\Delta V\Vert_{\mathcal{Z}^{2}_{\mathbb{F}}(t,T)}^{2}\\
&\le18\varepsilon^{1-\alpha}\psi(K_{1})^{2}(T^{\alpha}+6+6\alpha c_{2}^{2}K_{2}+2K_{2}^{\alpha})\Vert\Delta V\Vert_{\mathcal{Z}^{2}_{\mathbb{F}}(t,T)}^{2}.
\end{aligned}
\end{equation*}
Considering (\ref{eq 3.24})-(\ref{eq 3.27}), we get that, for all $1\le i \le n$ and $t\in[T-\varepsilon,T]$,
\begin{equation*}\label{eq 3.28}\tag{3.28}
\begin{aligned}
&\hspace{2em}|\Delta Y_{t}^{i}|^{2}+\mathbb{E}_{t}^{\mathbb{Q}^{i}}\Big[\int_{t}^{T}|\Delta Z_{s}^{i}|^{2}ds\Big]\\
&\le 72\varepsilon\psi(K_{1})^{2}(T+6c_{2}^{2}K_{2}+2K_{2})\Vert\Delta U\Vert_{\mathcal{S}^{\infty}_{\mathbb{F}}(t,T)}^{2}+\widetilde{C}\varepsilon^{1-\alpha}\Vert\Delta V\Vert_{\mathcal{Z}^{2}_{\mathbb{F}}(t,T)}^{2},\\
\end{aligned}
\end{equation*}
where
\begin{equation}\label{eq 3.29}\tag{3.29}
\widetilde{C}=18n\psi(K_{1})^{2}L_{4}^{2}c_{2}^{2}(T^{\alpha}+8+6\alpha L_{4}^{2}c_{2}^{2}K_{2}+2\alpha K_{2})+18\psi(K_{1})^{2}(T^{\alpha}+6+6\alpha c_{2}^{2}K_{2}+2K_{2}^{\alpha}).
\end{equation}
In view of Lemma 2.3, we finally derive that
\begin{equation*}\label{eq 3.30}\tag{3.30}
\begin{aligned}
&\hspace{2em}\Vert\Delta Y\Vert_{\mathcal{S}^{\infty}_{\mathbb{F}}(T-\varepsilon,T)}^{2}+c_{1}^{2}\Vert\Delta Z\Vert_{\mathcal{Z}^{2}_{\mathbb{F}}(T-\varepsilon,T)}^{2}\\
&\le 144n\varepsilon\psi(K_{1})^{2}(T+6c_{2}^{2}K_{2}+2K_{2})\Vert\Delta U\Vert_{\mathcal{S}^{\infty}_{\mathbb{F}}(T-\varepsilon,T)}^{2}+2\widetilde{C}n\varepsilon^{1-\alpha}\Vert\Delta V\Vert_{\mathcal{Z}^{2}_{\mathbb{F}}(T-\varepsilon,T)}^{2}.
\end{aligned}
\end{equation*}
Therefore, taking $\varepsilon$ sufficiently small so that the mapping $\Psi$ is a contraction on the previously given set $\mathcal{B}_{\varepsilon}(K_{1},K_{2})$, the Banach fixed point theorem guarantees that the mapping $\Psi$ admits a unique fixed point. This concludes the proof.
\end{proof}

\section{Global solution with bounded terminal value}

\hspace{2.6em} In this section, our goal is to explore the solvability of global solutions for diagonally quadratic mean-field BSDE (\ref{eq 3.1}) with the help of the existence and the uniqueness result obtained for local solutions of mean-field BSDEs with diagonally quadratic generators.

Let $\xi:\Omega\rightarrow \mathbb{R}^{n}$ be $\mathcal{F}_{T}$-measurable and $f:\Omega\times[0,T]\times\mathbb{R}^{n}\times\mathbb{R}^{n\times d}\times\mathcal{P}_{2}(\mathbb{R}^{n+n\times d})\rightarrow \mathbb{R}^{n}$ be $\mathbb{F}$-progressively measurable. They satisfy the following three assumptions:

\noindent\textbf{(A4)} There exists $M_{1}\ge 0$ such that $\Vert \xi \Vert_{\infty}\le M_{1}$.

\noindent\textbf{(A5)} For $1\le i\le n,$ and $d\mathbb{P}\times dt$-almost all $(\omega,t)\in \Omega\times [0,T],$ there exist constants $L>0,\gamma>0,M_{3}>0$ and an $\mathbb{F}$-progressively measurable process $\zeta=(\zeta_{t})$ with $\big\Vert \int_{0}^{T}|\zeta_{t}(\omega)|^{2}dt\big\Vert_{\infty}\le M_{3}$, such that, for all $y\in \mathbb{R}^{n},$ $z\in \mathbb{R}^{n\times d},$ and $\mu\in \mathcal{P}_{2}(\mathbb{R}^{n+n\times d}),$
\begin{equation*}
|f^{i}(\omega,t,y,z,\mu)|\le \zeta_{t}(\omega)+L|y|+\frac{\gamma}{2}|z^{i}|^{2}+LW_{2}(\mu_{1},\delta_{0}).
\end{equation*}

\noindent\textbf{(A6)} For $1\le i\le n$, and $d\mathbb{P}\times dt$-almost all $(\omega,t)\in \Omega\times [0,T],$ we have that, for all $y,\bar{y}\in \mathbb{R}^{n},$ $z,\bar{z}\in \mathbb{R}^{n\times d},$ and $\mu,\bar{\mu}\in \mathcal{P}_{2}(\mathbb{R}^{n+n\times d}),$
\begin{equation*}
\begin{aligned}
&|f^{i}(\omega,t,y,z,\mu)-f^{i}(\omega,t,\bar{y},\bar{z},\bar{\mu})|\\
&\hspace{-2em}\le\ L\big(|y-\bar{y}|+W_{2}(\mu_{1},\bar{\mu}_{1})\big)+\psi\big(|y|\vee|\bar{y}|\vee W_{2}(\mu_{1},\delta_{0}) \vee W_{2}(\bar{\mu}_{1},\delta_{0})\big)\\
&\hspace{1em}\times\Big[\big(1+|z|+|\bar{z}|+W_{2}(\mu_{2},\delta_{0})+W_{2}(\bar{\mu}_{2},\delta_{0})\big)|z^{i}-\bar{z}^{i}|+\sum\limits_{j\ne i}|z^{j}-\bar{z}^{j}|+W_{2}(\mu_{2},\bar{\mu}_{2})\Big],\\[-1mm]
&|f^{i}(\omega,t,y,z,\mu)-f^{i}(\omega,t,y,z,\delta_{0})|\le L\big(1+W_{2}(\mu_{1},\delta_{0})\big).
\end{aligned}
\end{equation*}


\begin{remark}
In \cite{HH22-1}, the solvability of global solutions for one-dimensional mean-field BSDEs with quadratic growth was considered, while our assumptions (A5)-(A6) can be extended to the multi-dimensional case with diagonally quadratic generators. Moreover, provided that assumption (A5) holds, the generator $f$ is bounded with respect to $\mu_{2}=\mu(\mathbb{R}^{n}\times \cdot)\in \mathcal{P}_{2}(\mathbb{R}^{n\times d})$. For example, for $1\le i\le n$, the following generator
\begin{equation*}\label{eq 4.1}\tag{4.1}
f^{i}(\omega,t,y,z,\mu)=1+|y|+|z^{i}|^{2}+\sum_{j\neq i}\sin(|z^{j}|)+W_{2}(\mu_{1},\delta_{0})\cos\big(W_{2}(\mu_{2},\delta_{0})\big)
\end{equation*}
satisfies assumptions (A5)-(A6) for every $(t,y,z)\in [0,T]\times \mathbb{R}^{n}\times \mathbb{R}^{n\times d}$ and $\mu\in\mathcal{P}_{2}(\mathbb{R}^{n+n\times d})$.
\end{remark}

\begin{theorem}
Let the assumptions (A4)-(A6) hold. Then mean-field BSDE (\ref{eq 3.1}) has a unique global solution $(Y,Z)\in \mathcal{S}_{\mathbb{F}}^{\infty}(0,T;\mathbb{R}^{n})\times \mathcal{Z}_{\mathbb{F}}^{2}(0,T;\mathbb{R}^{n\times d})$ on the whole interval $[0,T]$. Moreover, there exist two positive constants $J_{1}$ and $J_{2}$ depending only on $(L,M_{1},M_{3},T)$ such that
\begin{equation*}\label{eq 4.2}\tag{4.2}
\Vert Y\Vert_{\mathcal{S}_{\mathbb{F}}^{\infty}(0,T)}\le J_{1}\ \ \text{and}\ \ \Vert Z\Vert_{\mathcal{Z}_{\mathbb{F}}^{2}(0,T)}^{2}\le J_{2}.
\end{equation*}
\end{theorem}

\begin{proof} We shall begin to analyze the solvability of BSDE (\ref{eq 3.1}) on some interval $[T-\delta_{\kappa},T]$, with a positive constant $\delta_{\kappa}$ that will be determined later. For this purpose, we define
\begin{equation*}\label{eq 4.3}\tag{4.3}
\widetilde{C}=M_{1}^{2}+M_{3}+3L^{2}+2.
\end{equation*}
Furthermore, let $\eta(\cdot)$ be the unique solution of the following ordinary differential equation:
\begin{equation*}\label{eq 4.4}\tag{4.4}
\frac{\eta(t)}{n}=\widetilde{C}+\widetilde{C}\int_{t}^{T}\big(2\eta(s)+1\big)ds+\widetilde{C}\int_{t}^{T}\frac{\eta(s)}{n}ds,\ \ t\in[0,T].
\end{equation*}
It is easy to observe that $\eta(\cdot)$ is a continuously decreasing function in $t$, and we define
\begin{equation*}
\kappa:=\sup\limits_{t\in[0,T]}\eta(t)=\eta(0).
\end{equation*}

Thanks to the constraint of $\Vert \xi\Vert_{\infty}^{2}\le n\widetilde{C}\le \kappa$, Theorem 3.1 shows that there exists a positive constant $\delta_{\kappa}$ which depends only on $\kappa$, such that BSDE (\ref{eq 3.1}) admits a unique local solution, $(Y,Z)$ on the time interval $[T-\delta_{\kappa},T]$, and it can be constructed through the Picard iteration.

Let us consider the Picard iteration: For $1\le i\le n$, we consider
\begin{equation*}
Y_{t}^{0,i}=\xi^{i}+\int_{t}^{T}Z_{s}^{0,i}dW_{s},\ \ t\in[T-\delta_{\kappa},T],
\end{equation*}
and for $j\ge 0$,
\begin{equation*}\label{eq 4.5}\tag{4.5}
\begin{aligned}
Y_{t}^{j+1,i}&=\xi^{i}+\int_{t}^{T}f^{i}(s,Y_{s}^{j},Z_{s}^{j}(Z_{s}^{j+1,i};i),\mathbb{P}_{(Y_{s}^{j},Z_{s}^{j})})ds-\int_{t}^{T}Z_{s}^{j+1,i}dW_{s},\\
&=\xi^{i}+\int_{t}^{T}(V_{1,s}^{j+1,i}+V_{2,s}^{j+1,i})ds-\int_{t}^{T}Z_{s}^{j+1,i}dW_{s},\ \ t\in[T-\delta_{\kappa},T],
\end{aligned}
\end{equation*}
where
\begin{equation*}
V_{1,s}^{j+1,i}=f^{i}(s,0,Z_{s}^{j}(0;i),\delta_{0})+f^{i}(s,Y_{s}^{j},Z_{s}^{j}(Z_{s}^{j+1,i};i),\mathbb{P}_{(Y_{s}^{j},Z_{s}^{j})})-f^{i}(s,0,Z_{s}^{j}(Z_{s}^{j+1,i};i),\delta_{0})
\end{equation*}
and
\begin{equation*}
V_{2,s}^{j+1,i}=f^{i}(s,0,Z_{s}^{j}(Z_{s}^{j+1,i};i),\delta_{0})-f^{i}(s,0,Z_{s}^{j}(0;i),\delta_{0}).
\end{equation*}
By virtue of the assumption (A6), we can deduce the existence of a process $\Pi_{s}^{j+1,i}$ such that
$$V_{2,s}^{j+1,i}=Z_{s}^{j+1,i}\Pi_{s}^{j+1,i}\ \ \text{and}\ \ |\Pi_{s}^{j+1,i}|\le \psi(0)(1+2|Z_{s}^{j}|+|Z_{s}^{j+1}|).$$
Then the process
\begin{equation*}
\widetilde{W}_{t}^{j+1,i}:=W_{t}-\int_{0}^{t}\Pi_{s}^{j+1,i}{\bm{1}}_{[T-\delta_{\kappa},T]}(s)ds,\ \ t\in[0,T]
\end{equation*}
is a Brownian motion corresponding to an equivalent probability measure $\widetilde{\mathbb{Q}}^{j+1,i}$ defined by
\begin{equation*}
d\widetilde{\mathbb{Q}}^{j+1,i}:=\mathscr{E}(\Pi^{j+1,i}{\bm{1}}_{[T-\delta_{\kappa},T]}\cdot W)_{0}^{T}d\mathbb{P},
\end{equation*}
which is denoted by $\widetilde{\mathbb{Q}}^{i}$ hereafter for simplicity, and whose expectation is denoted by $\widetilde{\mathbb{E}}^{i}$ (note that $\mathbb{E}^{\mathbb{P}}[\mathscr{E}(\Pi^{j+1,i}{\bm{1}}_{[T-\delta_{\kappa},T]}\cdot W)_{0}^{T}]=1$, refer to \cite[Theorem 2.3]{K94}). Therefore, we can rewrite BSDE (\ref{eq 4.5}): For $1\le i\le n$ and $j\ge0$,
\begin{equation*}
Y_{t}^{j+1,i}=\xi^{i}+\int_{t}^{T}V_{1,s}^{j+1,i}ds-\int_{t}^{T}Z_{s}^{j+1,i}d\widetilde{W}_{s}^{j+1,i}.
\end{equation*}

By induction, we will show the following inequality: For $1\le i\le n$, $j\ge0$,
\begin{equation*}\label{eq 4.6}\tag{4.6}
|Y_{t}^{j,i}|^{2}\le \frac{\eta(t)}{n},\ \ t\in[T-\delta_{\kappa},T].
\end{equation*}
Actually, it's obvious that $|Y_{t}^{0,i}|^{2}\le \frac{\eta(t)}{n}$, and assume that $|Y_{t}^{j,i}|^{2}\le \frac{\eta(t)}{n}$, for $t\in [T-\delta_{\kappa},T]$. Then we just need to prove that $|Y_{t}^{j+1,i}|^{2}\le \frac{\eta(t)}{n}$, for $t\in [T-\delta_{\kappa},T]$. Applying It\^{o}'s formula to $|Y_{t}^{j+1,i}|^{2}$ and combining the assumptions (A4)-(A6) and Fubini's theorem, we derive that, for $r\in[T-\delta_{\kappa},t]$,
\begin{equation*}
\begin{aligned}
&\hspace{1.3em}\widetilde{\mathbb{E}}^{i}_{r}[|Y_{t}^{j+1,i}|^{2}]+\widetilde{\mathbb{E}}^{i}_{r}\Big[\int_{t}^{T}|Z_{s}^{j+1,i}|^{2}ds\Big]\\
&\le \widetilde{\mathbb{E}}^{i}_{r}[|\xi^{i}|^{2}]+\widetilde{\mathbb{E}}^{i}_{r}\Big[\int_{t}^{T}\Big(2|Y_{s}^{j+1,i}|^{2}+|\zeta_{s}(\omega)|^{2}+3L^{2}\big(|Y_{s}^{j}|^{2}+\Vert Y_{s}^{j}\Vert_{L^{2}(\Omega)}^{2}+1\big)\Big)ds\Big]\\
&\le M_{1}^{2}+M_{3}+3L^{2}\int_{t}^{T}\Big(\widetilde{\mathbb{E}}^{i}_{r}[|Y_{s}^{j}|^{2}]+\Vert Y_{s}^{j}\Vert_{L^{2}(\Omega)}^{2}+1\Big)ds+2\int_{t}^{T}\widetilde{\mathbb{E}}^{i}_{r}[|Y_{s}^{j+1,i}|^{2}]ds.
\end{aligned}
\end{equation*}
This together with the definition of $\widetilde{C}$ and the inequality (\ref{eq 4.6}) for $Y^{j}$ yields that
\begin{equation*}
\widetilde{\mathbb{E}}^{i}_{r}[|Y_{t}^{j+1,i}|^{2}]\le \widetilde{C}+\widetilde{C}\int_{t}^{T}\big(2\eta(s)+1\big)ds+\widetilde{C}\int_{t}^{T}\widetilde{\mathbb{E}}^{i}_{r}[|Y_{s}^{j+1,i}|^{2}]ds.
\end{equation*}
It follows from (\ref{eq 4.4}) and the comparison theorem that
\begin{equation*}
\widetilde{\mathbb{E}}^{i}_{r}[|Y_{t}^{j+1,i}|^{2}]\le \frac{\eta(t)}{n},\ \ t\in[T-\delta_{\kappa},T].
\end{equation*}
By setting $r=t$, we obtain
\begin{equation*}
|Y_{t}^{j+1,i}|^{2}\le \frac{\eta(t)}{n},\ \ t\in[T-\delta_{\kappa},T].
\end{equation*}
Therefore, inequality (\ref{eq 4.6}) holds. Since $Y_{t}^{i}=\lim\limits_{j\rightarrow \infty}Y_{t}^{j,i}$, our constructed local solution $(Y,Z)$ in $[T-\delta_{\kappa},T]$ satisfies the following estimate:
\begin{equation*}
|Y_{t}^{i}|^{2}\le \frac{\eta(t)}{n}\ \ \text{and}\ \ |Y_{t}|^{2}\le \eta(t),\ \ t\in[T-\delta_{\kappa},T].
\end{equation*}
Especially, $|Y_{T-\delta_{\kappa}}|^{2}\le \eta(T-\delta_{\kappa})\le \kappa$.

Taking $T-\delta_{\kappa}$ as the terminal time and $Y_{T-\delta_{\kappa}}$ as the terminal value and applying Theorem 3.1, we obtain that mean-field BSDE (\ref{eq 3.1}) has a local solution $(Y,Z)$ on $[T-2\delta_{\kappa},T-\delta_{\kappa}]$ through the Picard iteration. Also, using the Picard iteration and the fact that $|Y_{T-\delta_{\kappa}}^{i}|^{2}\le \frac{\eta(T-\delta_{\kappa})}{n}$, we derive that
\begin{equation*}
|Y_{t}^{i}|^{2}\le \frac{\eta(t)}{n},\ \ t\in[T-2\delta_{\kappa},T-\delta_{\kappa}].
\end{equation*}
Repeating the preceding process, we can extend the pair $(Y,Z)$ to the whole interval $[0,T]$ within a finite number of steps such that $Y$ is uniformly bounded by $\kappa$. Namely, there is a positive constant $J_{1}$ which depends only on $(L,M_{1},M_{3},T)$ such that
\begin{equation*}\label{eq 4.7}\tag{4.7}
\Vert Y\Vert_{\mathcal{S}_{\mathbb{F}}^{\infty}(0,T)}\le J_{1}.
\end{equation*}

We now show that $Z\cdot W$ is a $BMO(\mathbb{P})$-martingale. First, we define the following function:
\begin{equation*}\label{eq 4.8}\tag{4.8}
\phi(x):=\frac{1}{\gamma^{2}}\big(\exp\{\gamma |x|\}-\gamma |x|-1\big),\ \ x\in \mathbb{R}.
\end{equation*}
Then it is easy to see that, for $x\in \mathbb{R}$,
\begin{equation*}\label{eq 4.9}\tag{4.9}
\phi^{\prime}(x)=\frac{1}{\gamma}\big(\exp\{\gamma |x|\}-1\big)\text{sgn}(x),\ \ \phi^{\prime \prime}(x)=\exp\{\gamma |x|\},\ \ \phi^{\prime \prime}(x)-\gamma|\phi^{\prime}(x)|=1.
\end{equation*}
Using It\^{o}'s formula to compute $\phi(Y_{t}^{i})$ and combining this with assumptions (A5) and (\ref{eq 4.9}), we deduce that
\begin{equation*}
\begin{aligned}
&\hspace{2em}\phi(Y_{t}^{i})+\frac{1}{2}\mathbb{E}_{t}\Big[\int_{t}^{T}|Z_{s}^{i}|^{2}ds\Big]\\
&\le \phi(\Vert\xi^{i}\Vert_{\infty})+\mathbb{E}_{t}\Big[\int_{t}^{T}|\phi^{\prime}(Y_{s}^{i})|\big(\zeta_{s}(\omega)+L|Y_{s}|+L\Vert Y_{s}\Vert_{L^{2}(\Omega)}\big)ds\Big]\\
&\le \phi(M_{1})+\phi^{\prime}(J_{1})\big(\sqrt{TM_{3}}+2LJ_{1}T\big).
\end{aligned}
\end{equation*}
Consequently, we have
\begin{equation*}
\Vert Z\Vert_{\mathcal{Z}^{2}_{\mathbb{F}}(0,T)}^{2}=\Vert Z\cdot W\Vert_{BMO(\mathbb{P})}^{2}\le 2n\phi(M_{1})+2n\phi^{\prime}(J_{1})\big(\sqrt{TM_{3}}+2LJ_{1}T\big).
\end{equation*}
This together with (\ref{eq 4.7}) implies (\ref{eq 4.2}).

Finally, we shall prove the uniqueness. Let $(Y,Z)$ and $(\widetilde{Y},\widetilde{Z})$ be two adapted solutions of mean-field BSDE (\ref{eq 3.1}). We put
\begin{equation*}
\Delta Y:=Y-\widetilde{Y},\ \ \Delta Z:=Z-\widetilde{Z}.
\end{equation*}
Then for $1\le i\le n$, we have that, for arbitrary $\varepsilon>0$ and $t\in [T-\varepsilon,T]$,
\begin{equation*}
\Delta Y_{t}^{i}+\int_{t}^{T}\Delta Z_{s}^{i}d\widetilde{W}_{s}^{i}=\int_{t}^{T}\Big(f^{i}(s,Y_{s},Z_{s}(\widetilde{Z}_{s}^{i};i),\mathbb{P}_{(Y_{s},Z_{s})})-f^{i}(s,\widetilde{Y}_{s},\widetilde{Z}_{s},\mathbb{P}_{(\widetilde{Y}_{s},\widetilde{Z}_{s})})\Big)ds,
\end{equation*}
where $\widetilde{W}^{i}$ is defined analogously to that of (\ref{eq 3.21}). Similar to (\ref{eq 3.22})-(\ref{eq 3.30}), we conclude that, for all $t\in [T-\varepsilon,T]$,
\begin{equation*}
\begin{aligned}
&\hspace{2em}\Vert\Delta Y\Vert_{\mathcal{S}^{\infty}_{\mathbb{F}}(t,T)}^{2}+c_{1}^{2}\Vert\Delta Z\Vert_{\mathcal{Z}^{2}_{\mathbb{F}}(t,T)}^{2}\\
&\le 24nL^{2}\varepsilon\Vert\Delta Y\Vert_{\mathcal{S}^{\infty}_{\mathbb{F}}(t,T)}^{2}+6n\psi(J_{1})^{2}(nc_{2}^{2}+1)\varepsilon\Vert\Delta Z\Vert_{\mathcal{Z}^{2}_{\mathbb{F}}(t,T)}^{2}.
\end{aligned}
\end{equation*}
Thus, if $\varepsilon$ is small enough, we get $Y=\widetilde{Y}$ and $Z=\widetilde{Z}$ on the time interval $[T-\varepsilon,T]$ immediately. Repeating iteratively within a finite number of times, the uniqueness on the whole time interval $[0,T]$ can be obtained. The proof is then complete.
\end{proof}
\section{Global solution with unbounded terminal value}
\hspace{2.6em} We now investigate the existence and the uniqueness of solutions for diagonally quadratic mean-field BSDEs with unbounded terminal values and convex generators. The form of the equation is expressed as follows:
\begin{equation}\label{eq 5.1}\tag{5.1}
	Y_{t}=\xi+\int_{t}^{T}f(s,Y_{s},Z_{s},\mathbb{P}_{Y_{s}})ds-\int_{t}^{T}Z_{s}dW_s,\ t\in[0,T].
\end{equation}

Let $\xi:\Omega\rightarrow \mathbb{R}^{n}$ be $\mathcal{F}_{T}$-measurable, and let $f:\Omega\times[0,T]\times\mathbb{R}^{n}\times\mathbb{R}^{n\times d}\times\mathcal{P}_{1}(\mathbb{R}^{n})\rightarrow \mathbb{R}^{n}$ be $\mathbb{F}$-progressively measurable, satisfying the following assumptions:

\noindent\textbf{(A7)} $\xi$ has exponential moments of all orders, i.e., for all $p\ge1$, we have $\mathbb{E}\big[\exp\big\{p|\xi|\big\}\big]<+\infty.$

\noindent\textbf{(A8)} For $1\le i\le n,\ f^{i}(\omega,t,y,z,\mu)$ varies with $(\omega,t,y,z^{i},\mu)$ only, and grows linearly in $(y,\mu)$ and quadratically in $z^{i}$, i.e., $d\mathbb{P}\times dt\text{-}$almost all $(\omega,t)\in \Omega\times [0,T]$, there exist constants $K>0$ and $\gamma>0$, such that, for all $y\in\mathbb{R}^{n}$, $z\in\mathbb{R}^{n\times d}$, and $\mu\in\mathcal{P}_{1}(\mathbb{R}^{n})$,
\begin{equation*}
|f^{i}(\omega,t,y,z,\mu)|\le \zeta_{t}(\omega)+K|y|+\frac{\gamma}{2}|z^{i}|^{2}+K W_{1}(\mu,\delta_{0}).
\end{equation*}
Here $\zeta$ is an $\mathbb{F}$-progressively measurable process s.t., for all $p\ge1$, $\mathbb{E}\big[\exp\big\{p\int_{0}^{T}\zeta_{t}(\omega)dt\big\}\big]<+\infty.$

\noindent\textbf{(A9)} The function $f(\omega,t,y,z,\mu)$ is uniformly Lipschitz continuous with respect to $(y,\mu)$, i.e., $d\mathbb{P}\times dt\text{-}$almost all $(\omega,t)\in \Omega\times [0,T]$, we have that, for all $y,\bar{y}\in \mathbb{R}^{n},$ $z\in \mathbb{R}^{n\times d},$ and $\mu,\bar{\mu}\in \mathcal{P}_{1}(\mathbb{R}^{n}),$
\begin{equation*}
|f(\omega,t,y,z,\mu)-f(\omega,t,\bar{y},z,\bar{\mu})|\le K\big(|y-\bar{y}|+W_{1}(\mu,\bar{\mu})\big).
\end{equation*}

\noindent\textbf{(A10)} For all $1\le i\le n$, and $(y,\mu)\in \mathbb{R}^{n}\times\mathcal{P}_{1}(\mathbb{R}^{n})$, it holds that $f^{i}(\omega,t,y,\cdot,\mu)$ is either convex or concave.

\begin{remark} The assumption (A10) is valid for the generator $f$, provided that certain components of $f$ are convex in $z$ while others are concave in $z$.
\end{remark}
To begin with, we present two lemmas which provide some bounds on the (possibly unbounded) solutions of scalar-valued quadratic BSDEs. For the detailed proof, the reader can refer to Fan, Hu and Tang \cite[Proposition 1]{FH20}.
\begin{lemma} Assume that the generator $f:\Omega\times [0,T]\times \mathbb{R}^{1\times d}\rightarrow \mathbb{R}$ is $\mathbb{F}$-progressively measurable and there exists an $\mathbb{F}$-progressively measurable non-negative process $(\bar{\zeta}_{t})_{t\in[0,T]}$ such that, $d\mathbb{P}\times dt\text{-}a.e.$, $(\omega,t)\in \Omega\times [0,T]$, for all $z\in\mathbb{R}^{1\times d}$,
$$|f(\omega,t,z)|\le \bar{\zeta}_{t}(\omega)+\frac{\gamma}{2}|z|^{2}.$$
Then, for any solution $(Y,Z)$ of BSDE (\ref{eq 3.2}) satisfying
\begin{equation*}
\mathbb{E}\bigg[\exp\bigg\{2\gamma\sup\limits_{t\in[0,T]}|Y_{t}|+2\gamma\int_{0}^{T}\bar{\zeta}_{s}(\omega)ds\bigg\}\bigg]< +\infty,
\end{equation*}
we have that for all $t\in [0,T]$,
\begin{equation*}
\exp\big\{\gamma |Y_{t}|\big\}\le \mathbb{E}_{t}\bigg[\exp\bigg\{\gamma |\xi|+\gamma\int_{t}^{T}\bar{\zeta}_{s}(\omega)ds\bigg\}\bigg].
\end{equation*}
\end{lemma}

\begin{lemma} Assume that the generator $f:\Omega\times [0,T]\times \mathbb{R}^{1\times d}\rightarrow \mathbb{R}$ is $\mathbb{F}$-progressively measurable and there exists an $\mathbb{F}$-progressively measurable non-negative process $(\bar{\zeta}_{t})_{t\in[0,T]}$ such that, $d\mathbb{P}\times dt\text{-}a.e.$, $(\omega,t)\in \Omega\times [0,T]$, for all $z\in\mathbb{R}^{1\times d}$,
$$f(\omega,t,z)\le \bar{\zeta}_{t}(\omega)+\frac{\gamma}{2}|z|^{2}.$$
Then, for any solution $(Y,Z)$ of BSDE (\ref{eq 3.2}) satisfying
\begin{equation*}
\mathbb{E}\bigg[\exp\bigg\{2\gamma\sup\limits_{t\in[0,T]}Y_{t}^{+}+2\gamma\int_{0}^{T}\bar{\zeta}_{s}(\omega)ds\bigg\}\bigg]< +\infty,
\end{equation*}
we have that for all $t\in [0,T]$,
\begin{equation*}
\exp\big\{\gamma Y_{t}^{+}\big\}\le \mathbb{E}_{t}\bigg[\exp\bigg\{\gamma \xi^{+}+\gamma\int_{t}^{T}\bar{\zeta}_{s}(\omega)ds\bigg\}\bigg].
\end{equation*}
\end{lemma}

\begin{theorem} Let the assumptions (A7)-(A10) hold. Then mean-field BSDE (\ref{eq 5.1}) possesses a unique global solution $(Y,Z)\in \mathcal{E}(0,T;\mathbb{R}^{n})\times \mathcal{M}(0,T;\mathbb{R}^{n\times d})$ on the whole interval $[0,T]$.
\end{theorem}

\begin{proof} \textbf{Step 1.}\ \ For a pair of processes $(U,V)\in \mathcal{E}(0,T;\mathbb{R}^{n})\times \mathcal{M}(0,T;\mathbb{R}^{n\times d})$, we consider the following decoupled system of quadratic BSDEs:
\begin{equation}\label{eq 5.2}\tag{5.2}
Y_{t}^{i}=\xi^{i}+\int_{t}^{T}f^{i}(s,U_{s},Z_{s}^{i},\mathbb{P}_{U_{s}})ds-\int_{t}^{T}Z_{s}^{i}dW_{s},\ \ t\in[0,T],\ \ 1\le i\le n.
\end{equation}
In view of assumption (A8), we can easily get that, $d\mathbb{P}\times dt\text{-}a.e.$, for $z\in \mathbb{R}^{1\times d}$,
\begin{equation*}\label{eq 5.3}\tag{5.3}
|f^{i}(s,U_{s},z,\mathbb{P}_{U_{s}})|\le \zeta_{s}(\omega)+K|U_{s}|+K\Vert U_{s}\Vert_{L^{1}(\Omega)}+\frac{\gamma}{2}|z|^{2}.
\end{equation*}
By applying H\"{o}lder's inequality and Jensen's inequality, and taking into account assumption (A7) as well as the fact that $U\in \mathcal{E}(0,T;\mathbb{R}^{n})$, we deduce that, for $1\le i\le n$,
\begin{equation*}
\mathbb{E}\bigg[\exp\bigg\{q\Big(|\xi^{i}|+\int_{0}^{T}\big(\zeta_{s}(\omega)+K|U_{s}|+K\Vert U_{s}\Vert_{L^{1}(\Omega)}\big)ds\Big)\bigg\}\bigg]<+\infty,\ \  q>1.
\end{equation*}
Combining with the assumption (A10), then it can be inferred from Corollary 6 as presented in the work of Briand and Hu \cite{BH08} that the system of BSDEs (\ref{eq 5.2}) admits a unique solution $(Y,Z)\in \mathcal{E}(0,T;\mathbb{R}^{n})\times \mathcal{M}(0,T;\mathbb{R}^{n\times d})$.

\noindent\textbf{Step 2.}\ \ Thus, we can put $(Y^{0},Z^{0})=(0,0)$ and recursively define the sequence of processes $\{(Y^{m},Z^{m})\}_{m=1}^{\infty}$ in the space of processes $\mathcal{E}(0,T;\mathbb{R}^{n})\times \mathcal{M}(0,T;\mathbb{R}^{n\times d})$ by the unique solution of the system of quadratic BSDEs
\begin{equation}\label{eq 5.4}\tag{5.4}
Y_{t}^{m+1;i}=\xi^{i}+\int_{t}^{T}f^{i}(s,Y_{s}^{m},Z_{s}^{m+1;i},\mathbb{P}_{Y_{s}^{m}})ds-\int_{t}^{T}Z_{s}^{m+1;i}dW_{s},\ \ t\in[0,T],\ \ 1\le i\le n,
\end{equation}
where for the sake of convenience, we shall denote the $i$th component of $Y^{m}$ and the $i$th row of $Z^{m}$ by $Y^{m;i}$ and $Z^{m;i}$, respectively. Next, we are going to prove that $\{(Y^{m},Z^{m})\}_{m=1}^{\infty}$ is a Cauchy sequence in $\mathcal{S}^{q}_{\mathbb{F}}(0,T;\mathbb{R}^{n})\times \mathcal{H}^{q}_{\mathbb{F}}(0,T;\mathbb{R}^{n\times d})$, for every $q\ge1$ and thus converges to a pair of adapted processes $(Y,Z)\in\mathcal{E}(0,T;\mathbb{R}^{n})\times \mathcal{M}(0,T;\mathbb{R}^{n\times d})$, which is the unique solution of BSDE (\ref{eq 5.1}).

The main objective of step 2 is to show that, for all $q>1$,
\begin{equation}\label{eq 5.5}\tag{5.5}
\sup\limits_{m\ge0}\mathbb{E}\bigg[\exp\Big\{q\gamma\sup\limits_{t\in[0,T]}|Y_{t}^{m}|\Big\}\bigg]\le C(q),
\end{equation}
where
\begin{equation*}
\begin{aligned}
&C(q):=\big(R(2q)R(8nq)\big)^{[4nK T]+1}\mathbb{E}\Big[\exp\Big\{4n(8n)^{[4nK T]+1}q\gamma|\xi|\Big\}\Big]\\
&\hspace{4em}\cdot\mathbb{E}\bigg[\exp\bigg\{4n(16n)^{[4nK T]+1}q\gamma\int_{0}^{T}\zeta_{s}(\omega)ds\bigg\}\bigg]<+\infty
\end{aligned}
\end{equation*}
with
\begin{equation}\label{eq 5.6}\tag{5.6}
R(q):=\Big(\frac{q}{q-1}\Big)^{2q}.
\end{equation}
Indeed, from assumption (A8) we have that, for $1\le i\le n$ and $m\ge0$, $d\mathbb{P}\times dt\text{-}a.e.,$
\begin{equation*}\label{eq 5.7}\tag{5.7}
|f^{i}(s,Y_{s}^{m},z,\mathbb{P}_{Y_{s}^{m}})|\le \zeta_{s}(\omega)+K |Y_{s}^{m}|+K \Vert Y_{s}^{m}\Vert_{L^{1}(\Omega)}+\frac{\gamma}{2}|z|^{2}.
\end{equation*}
In view of the fact that $Y^{m}\in \mathcal{E}(0,T;\mathbb{R}^{n})$, by using H\"{o}lder's and Jensen's inequalities, we derive that, for all $q>1$,
\begin{equation*}
\mathbb{E}\bigg[\exp\bigg\{q\bigg(\sup\limits_{t\in[0,T]}|Y_{t}^{m+1;i}|+\int_{0}^{T}\Big(\zeta_{s}(\omega)+K |Y_{s}^{m}|+K \Vert Y_{s}^{m}\Vert_{L^{1}(\Omega)}\Big)ds\bigg)\bigg\}\bigg]<+\infty.
\end{equation*}
Subsequently, we can apply Lemma 5.1 to obtain that, for $1\le i\le n$ and $m\ge0$,
\begin{equation*}
\exp\big\{\gamma |Y_{t}^{m+1;i}|\big\}\le \mathbb{E}_{t}\bigg[\exp\bigg\{\gamma |\xi^{i}|+\gamma\int_{t}^{T}\Big(\zeta_{s}(\omega)+K |Y_{s}^{m}|+K \Vert Y_{s}^{m}\Vert_{L^{1}(\Omega)}\Big)ds\bigg\}\bigg],\ \ t\in [0,T].
\end{equation*}
From here we see immediately that, for $m\ge0$,
\begin{equation}\label{eq 5.8}\tag{5.8}
\exp\big\{\gamma |Y_{t}^{m+1}|\big\}\le \mathbb{E}_{t}\bigg[\exp\bigg\{n\gamma |\xi|+n\gamma\int_{t}^{T}\Big(\zeta_{s}(\omega)+K |Y_{s}^{m}|+K \Vert Y_{s}^{m}\Vert_{L^{1}(\Omega)}\Big)ds\bigg\}\bigg],\ \ t\in [0,T].
\end{equation}
By using Doob's maximal inequality, and H\"{o}lder's and Jensen's inequalities, we conclude that, for all $q>1$, $m\ge0$ and $t\in[0,T]$,
\begin{equation}\label{eq 5.9}\tag{5.9}
\begin{aligned}
&\hspace{2em}\mathbb{E}\bigg[\exp\bigg\{q\gamma \sup\limits_{s\in[t,T]}|Y_{s}^{m+1}|\bigg\}\bigg]\\
&\le \Big(\frac{q}{q-1}\Big)^{q}\mathbb{E}\bigg[\exp\bigg\{nq\gamma |\xi|+nq\gamma\int_{t}^{T}\zeta_{s}(\omega)ds\\
&\hspace{7em}+nq\gamma K\big(\sup\limits_{s\in[t,T]}|Y_{s}^{m}|+\mathbb{E}[\sup\limits_{s\in[t,T]}|Y_{s}^{m}|]\big)(T-t)\bigg\}\bigg]\\
&\le \sqrt{I(q)}\bigg(\mathbb{E}\bigg[\exp\bigg\{4nq\gamma K\sup\limits_{s\in[t,T]}|Y_{s}^{m}|(T-t)\bigg\}\bigg]\bigg)^{\frac{1}{2}},
\end{aligned}
\end{equation}
where
\begin{equation*}\label{eq 5.10}\tag{5.10}
I(q):=R(q)\mathbb{E}\Big[\exp\Big\{2nq\gamma|\xi|+2nq\gamma\int_{0}^{T}\zeta_{s}(\omega)ds\Big\}\Big]
\end{equation*}
with $R(q)$ defined as specified in (\ref{eq 5.6}).

For the case of $K=0$, it is obvious from (\ref{eq 5.9}) that
\begin{equation}\label{eq 5.11}\tag{5.11}
\sup\limits_{m\ge0}\mathbb{E}\bigg[\exp\Big\{q\gamma\sup\limits_{t\in[0,T]}|Y_{t}^{m}|\Big\}\bigg]\le \sqrt{I(q)}.
\end{equation}
Otherwise, by setting $\varepsilon=\frac{1}{4nK}>0$, we let $m_{0}$ be the unique positive integer such that $T-m_{0}\varepsilon\le 0< T-(m_{0}-1)\varepsilon$, or in other words,
\begin{equation}\label{eq 5.12}\tag{5.12}
4nKT\le m_{0} <4nKT+1.
\end{equation}
If $m_{0}=1$, then for $t\in[0,T]$, we know $4nq\gamma K(T-t)\le 4nq\gamma KT\le 4nq\gamma K\varepsilon = q\gamma$, and it follows from (\ref{eq 5.9}) that, for all $q>1$, $m\ge0$ and $t\in[0,T]$,
\begin{equation*}
\mathbb{E}\bigg[\exp\bigg\{q\gamma \sup\limits_{s\in[t,T]}|Y_{s}^{m+1}|\bigg\}\bigg]\le\sqrt{I(q)}\bigg(\mathbb{E}\bigg[\exp\bigg\{q\gamma\sup\limits_{s\in[t,T]}|Y_{s}^{m}|\bigg\}\bigg]\bigg)^{\frac{1}{2}}.
\end{equation*}
By induction, we observe that
\begin{equation*}
\mathbb{E}\bigg[\exp\bigg\{q\gamma \sup\limits_{s\in[t,T]}|Y_{s}^{m+1}|\bigg\}\bigg]\le\sqrt{I(q)}^{1+\frac{1}{2}+\cdots+\frac{1}{2^{m}}}\bigg(\mathbb{E}\bigg[\exp\bigg\{q\gamma\sup\limits_{s\in[t,T]}|Y_{s}^{0}|\bigg\}\bigg]\bigg)^{\frac{1}{2^{m+1}}}\le I(q).
\end{equation*}
Therefore, a renewed application of H\"{o}lder's inequality implies that, for all $q>1$,
\begin{equation*}\label{eq 5.13}\tag{5.13}
\sup\limits_{m\ge0}\mathbb{E}\bigg[\exp\bigg\{q\gamma \sup\limits_{t\in[0,T]}|Y_{t}^{m}|\bigg\}\bigg]\le R(q)\mathbb{E}\big[\exp\big\{4nq\gamma|\xi|\big\}\big]\mathbb{E}\Big[\exp\Big\{4nq\gamma\int_{0}^{T}\zeta_{s}(\omega)ds\Big\}\Big].
\end{equation*}
If $m_{0}=2$, then for $t\in[T-\varepsilon,T]$, we know $4nq\gamma K(T-t)\le 4nq\gamma K\varepsilon=q\gamma$. By following the similar argument in the case of $m_{0}=1$, we deduce that, for all $q>1$,
\begin{equation*}\label{eq 5.14}\tag{5.14}
\sup\limits_{m\ge0}\mathbb{E}\bigg[\exp\bigg\{q\gamma \sup\limits_{t\in[T-\varepsilon,T]}|Y_{t}^{m}|\bigg\}\bigg]\le R(q)\mathbb{E}\big[\exp\big\{4nq\gamma|\xi|\big\}\big]\mathbb{E}\Big[\exp\Big\{4nq\gamma\int_{0}^{T}\zeta_{s}(\omega)ds\Big\}\Big],
\end{equation*}
which leads to
\begin{equation*}\label{eq 5.15}\tag{5.15}
\sup\limits_{m\ge0}\mathbb{E}\bigg[\exp\bigg\{q\gamma |Y_{T-\varepsilon}^{m}|\bigg\}\bigg]\le R(q)\mathbb{E}\big[\exp\big\{4nq\gamma|\xi|\big\}\big]\mathbb{E}\Big[\exp\Big\{4nq\gamma\int_{0}^{T}\zeta_{s}(\omega)ds\Big\}\Big].
\end{equation*}
We now consider the following system of quadratic BSDEs:
\begin{equation*}
Y_{t}^{m+1;i}=Y_{T-\varepsilon}^{m+1;i}+\int_{t}^{T-\varepsilon}f^{i}(s,Y_{s}^{m},Z_{s}^{m+1;i},\mathbb{P}_{Y_{s}^{m}})ds-\int_{t}^{T-\varepsilon}Z_{s}^{m+1;i}dW_{s},\ \ t\in[0,T-\varepsilon],\ \ 1\le i\le n.
\end{equation*}
Equation (\ref{eq 5.15}) allows to use a similar argument as that employed for (\ref{eq 5.13}), in order to prove that, for all $q>1$,
\begin{equation*}
\begin{aligned}
&\hspace{2em}\sup\limits_{m\ge0}\mathbb{E}\bigg[\exp\bigg\{q\gamma \sup\limits_{t\in[0,T-\varepsilon]}|Y_{t}^{m}|\bigg\}\bigg]\\
&\le R(q)\sup\limits_{m\ge 0}\mathbb{E}\big[\exp\big\{4nq\gamma|Y_{T-\varepsilon}^{m}|\big\}\big]\mathbb{E}\Big[\exp\Big\{4nq\gamma\int_{0}^{T}\zeta_{s}(\omega)ds\Big\}\Big]\\
&\le R(q)R(4nq)\mathbb{E}\big[\exp\big\{16n^{2}q\gamma|\xi|\big\}\big]\mathbb{E}\Big[\exp\Big\{32n^{2}q\gamma\int_{0}^{T}\zeta_{s}(\omega)ds\Big\}\Big].
\end{aligned}
\end{equation*}
Then by applying H\"{o}lder's inequality combined with (\ref{eq 5.14}), we obtain that
\begin{equation*}
\begin{aligned}
&\hspace{2em}\sup\limits_{m\ge0}\mathbb{E}\bigg[\exp\bigg\{q\gamma \sup\limits_{t\in[0,T]}|Y_{t}^{m}|\bigg\}\bigg]\\
&\le R(2q)R(8nq)\mathbb{E}\big[\exp\big\{32n^{2}q\gamma|\xi|\big\}\big]\mathbb{E}\Big[\exp\Big\{64n^{2}q\gamma\int_{0}^{T}\zeta_{s}(\omega)ds\Big\}\Big].
\end{aligned}
\end{equation*}
Therefore, assuming that $m_{0}$ satisfies (\ref{eq 5.12}), we can conclude that, for all $q>1$,
\begin{equation*}
\begin{aligned}
&\hspace{2em}\sup\limits_{m\ge0}\mathbb{E}\bigg[\exp\bigg\{q\gamma \sup\limits_{t\in[0,T]}|Y_{t}^{m}|\bigg\}\bigg]\\
&\le \big(R(2q)R(8nq)\big)^{m_{0}-1}\mathbb{E}\big[\exp\big\{4n(8n)^{m_{0}-1}q\gamma|\xi|\big\}\big]\mathbb{E}\Big[\exp\Big\{4n(16n)^{m_{0}-1}q\gamma\int_{0}^{T}\zeta_{s}(\omega)ds\Big\}\Big].
\end{aligned}
\end{equation*}
This together with (\ref{eq 5.11}) and (\ref{eq 5.12}) yields the estimate (\ref{eq 5.5}).

\noindent\textbf{Step 3.}\ \ We show that
\begin{equation}\label{eq 5.16}\tag{5.16}
\sup\limits_{m\ge0}\mathbb{E}\bigg[\Big(\int_{0}^{T}|Z_{s}^{m}|^{2}ds\Big)^{\frac{q}{2}}\bigg]<+\infty,\ \  q>1.
\end{equation}
In fact, applying It\^{o}-Tanaka's formula to the term $\exp\{2\gamma|Y_{t}^{m+1;i}|\}$ and using inequality (\ref{eq 5.7}), we deduce that, for $1\le i\le n$ and $m\ge0$,
\begin{equation*}
\begin{aligned}
&\hspace{2em}\exp\{2\gamma|Y_{0}^{m+1;i}|\}+\gamma^{2}\int_{0}^{T}\exp\{2\gamma|Y_{s}^{m+1;i}|\}|Z_{s}^{m+1;i}|^{2}ds\\
&\le \exp\{2\gamma |\xi^{i}|\}+2\gamma\int_{0}^{T}\exp\{2\gamma|Y_{s}^{m+1;i}|\}\big(\zeta_{s}(\omega)+K |Y_{s}^{m}|+K\mathbb{E}[|Y_{s}^{m}|]\big)ds\\
&\hspace{1em}-2\gamma\int_{0}^{T}\exp\{2\gamma|Y_{s}^{m+1;i}|\}\text{sgn}(Y_{s}^{m+1;i})Z_{s}^{m+1;i}dW_{s}.
\end{aligned}
\end{equation*}
Then by using Burkholder-Davis-Gundy inequality, as well as Young's and Jensen's inequalities, we derive that, for all $q>1$,
\begin{equation*}\label{eq 5.17}\tag{5.17}
\begin{aligned}
&\hspace{2em}\mathbb{E}\bigg[\bigg(\int_{0}^{T}|Z_{s}^{m+1;i}|^{2}ds\bigg)^{\frac{q}{2}}\bigg]\\
&\le M\bigg(\mathbb{E}\bigg[\exp\bigg\{2q\gamma\sup\limits_{t\in[0,T]}|Y_{t}^{m+1;i}|\bigg\}\bigg]+\mathbb{E}\bigg[\Big(\int_{0}^{T}\big(\zeta_{s}(\omega)+K |Y_{s}^{m}|+K\mathbb{E}[|Y_{s}^{m}|]\big)ds\Big)^{q}\bigg]\bigg)\\
&\le \widetilde{M}\bigg(\sup\limits_{m\ge0}\mathbb{E}\bigg[\exp\Big\{2q\gamma\sup\limits_{t\in[0,T]}|Y_{t}^{m}|\Big\}\bigg]+\mathbb{E}\bigg[\exp\bigg\{\int_{0}^{T}\zeta_{s}(\omega)ds\bigg\}\bigg]\bigg),
\end{aligned}
\end{equation*}
where $M$ and $\widetilde{M}$ represent two constants which depend only on $(q,\gamma)$ and $(q,\gamma,K,T)$, respectively. Combining this estimate with (\ref{eq 5.5}), we conclude that (\ref{eq 5.16}) is satisfied.

\noindent\textbf{Step 4.}\ \ In what follows, without loss of generality, we suppose that the generator $f$ is component-wise convex in (A10), i.e., for all $1\le i\le n$ and $(y,\mu)\in \mathbb{R}^{n}\times \mathcal{P}_{1}(\mathbb{R}^{n})$, we have that, $d\mathbb{P}\times dt\text{-}a.e.$, $f^{i}(\omega,t,y,\cdot,\mu)$ is convex.

For any given $m,p\ge1$ and $\theta\in(0,1)$, let us define
\begin{equation*}
\Delta_{\theta}Y^{m,p}:=\frac{Y^{m+p}-\theta Y^{m}}{1-\theta}\ \ \text{and}\ \ \Delta_{\theta}Z^{m,p}:=\frac{Z^{m+p}-\theta Z^{m}}{1-\theta}.
\end{equation*}
The main objective of this step is to give an exponential estimate of $\Delta_{\theta}Y^{m,p}$. We observe that the pair of processes $(\Delta_{\theta}Y^{m,p},\Delta_{\theta}Z^{m,p})\in \mathcal{E}(0,T;\mathbb{R}^{n})\times \mathcal{M}(0,T;\mathbb{R}^{n\times d})$ satisfies the following system of quadratic BSDEs: For $1\le i\le n$,
\begin{equation}\label{eq 5.18}\tag{5.18}
\Delta_{\theta}Y_{t}^{m,p;i}=\xi^{i}+\int_{t}^{T}\Delta_{\theta}f^{m,p;i}(s,\Delta_{\theta}Z_{s}^{m,p;i})ds-\int_{t}^{T}\Delta_{\theta}Z_{s}^{m,p;i}dW_{s},\ \ t\in [0,T],
\end{equation}
where
\begin{equation*}\label{eq 5.19}\tag{5.19}
\Delta_{\theta}f^{m,p;i}(s,z):=\frac{1}{1-\theta}\Big(f^{i}(s,Y_{s}^{m+p-1},(1-\theta)z+\theta Z_{s}^{m;i},\mathbb{P}_{Y_{s}^{m+p-1}})-\theta f^{i}(s,Y_{s}^{m-1},Z_{s}^{m;i},\mathbb{P}_{Y_{s}^{m-1}})\Big).
\end{equation*}
Recalling the assumptions (A8)-(A10), we see that, $d\mathbb{P}\times dt\text{-}a.e.$, for all $z\in \mathbb{R}^{1\times d}$,
\begin{equation*}
\Delta_{\theta}f^{m,p;i}(s,z)\le \zeta_{s}(\omega)+K\big(|\Delta_{\theta}Y_{s}^{m-1,p}|+\mathbb{E}[|\Delta_{\theta}Y_{s}^{m-1,p}|]\big)+2K(|Y_{s}^{m-1}|+\mathbb{E}[|Y_{s}^{m-1}|])+\frac{\gamma}{2}|z|^{2},
\end{equation*}
which together with (\ref{eq 5.5}) and Jensen's inequality shows that all the conditions in Lemma 5.2 are satisfied for BSDE (\ref{eq 5.18}), and so we conclude that, for all $1\le i\le n$ and $t\in[0,T]$,
\begin{equation*}\label{eq 5.20}\tag{5.20}
\begin{aligned}
&\hspace{2em}\exp\Big\{\gamma\big(\Delta_{\theta}Y_{t}^{m,p;i}\big)^{+}\Big\}\\
&\le \mathbb{E}_{t}\bigg[\exp\bigg\{\gamma (\xi^{i})^{+}+\gamma\int_{t}^{T}\Big(\zeta_{s}(\omega)+K\big(|\Delta_{\theta}Y_{s}^{m-1,p}|+\mathbb{E}[|\Delta_{\theta}Y_{s}^{m-1,p}|]\big)\\
&\hspace{6em}+2K(|Y_{s}^{m-1}|+\mathbb{E}[|Y_{s}^{m-1}|])\Big)ds\bigg\}\bigg].
\end{aligned}
\end{equation*}
Let us also define
\begin{equation*}
\Delta_{\theta}\widetilde{Y}^{m,p}:=\frac{Y^{m}-\theta Y^{m+p}}{1-\theta}\ \ \text{and}\ \ \Delta_{\theta}\widetilde{Z}^{m,p}:=\frac{Z^{m}-\theta Z^{m+p}}{1-\theta}.
\end{equation*}
Similar to above, we derive that, for $1\le i\le n$ and $t\in[0,T]$,
\begin{equation*}\label{eq 5.21}\tag{5.21}
\begin{aligned}
&\hspace{2em}\exp\Big\{\gamma\big(\Delta_{\theta}\widetilde{Y}_{t}^{m,p;i}\big)^{+}\Big\}\\
&\le \mathbb{E}_{t}\bigg[\exp\bigg\{\gamma (\xi^{i})^{+}+\gamma\int_{t}^{T}\Big(\zeta_{s}(\omega)+K\big(|\Delta_{\theta}\widetilde{Y}_{s}^{m-1,p}|+\mathbb{E}[|\Delta_{\theta}\widetilde{Y}_{s}^{m-1,p}|]\big)\\
&\hspace{6em}+2K(|Y_{s}^{m+p-1}|+\mathbb{E}[|Y_{s}^{m+p-1}|])\Big)ds\bigg\}\bigg].
\end{aligned}
\end{equation*}
Notice the fact that
\begin{equation*}
\big(\Delta_{\theta}Y_{t}^{m,p;i}\big)^{-}\le \big(\Delta_{\theta}\widetilde{Y}_{t}^{m,p;i}\big)^{+}+2|Y_{t}^{m+p}|\ \ \text{and}\ \ \big(\Delta_{\theta}\widetilde{Y}_{t}^{m,p;i}\big)^{-}\le \big(\Delta_{\theta}Y_{t}^{m,p;i}\big)^{+}+2|Y_{t}^{m}|.
\end{equation*}
By combining (\ref{eq 5.20}) and (\ref{eq 5.21}) and applying Jensen's inequality, we can obtain that, for all $1\le i\le n$ and $t\in[0,T]$,
\begin{equation*}
\begin{aligned}
&\hspace{1em}\exp\Big\{\gamma|\Delta_{\theta}Y_{t}^{m,p;i}|\Big\}=\exp\Big\{\gamma\big(\Delta_{\theta}Y_{t}^{m,p;i}\big)^{+}\Big\}\cdot\exp\Big\{\gamma\big(\Delta_{\theta}Y_{t}^{m,p;i}\big)^{-}\Big\}\\
&\le \mathbb{E}_{t}\bigg[\exp\bigg\{2\gamma |\xi|+2\gamma |Y_{t}^{m+p}|+2\gamma\int_{t}^{T}\Big(\zeta_{s}(\omega)+K\big(|\Delta_{\theta}Y_{s}^{m-1,p}|+\mathbb{E}[|\Delta_{\theta}Y_{s}^{m-1,p}|]+|\Delta_{\theta}\widetilde{Y}_{s}^{m-1,p}|\\
&\hspace{2.6em}+\mathbb{E}[|\Delta_{\theta}\widetilde{Y}_{s}^{m-1,p}|]\big)+2K\big(|Y_{s}^{m+p-1}|+\mathbb{E}[|Y_{s}^{m+p-1}|]+|Y_{s}^{m-1}|+\mathbb{E}[|Y_{s}^{m-1}|]\big)\Big)ds\bigg\}\bigg]
\end{aligned}
\end{equation*}
and
\begin{equation*}
\begin{aligned}
&\hspace{1em}\exp\Big\{\gamma|\Delta_{\theta}\widetilde{Y}_{t}^{m,p;i}|\Big\}=\exp\Big\{\gamma\big(\Delta_{\theta}\widetilde{Y}_{t}^{m,p;i}\big)^{+}\Big\}\cdot\exp\Big\{\gamma\big(\Delta_{\theta}\widetilde{Y}_{t}^{m,p;i}\big)^{-}\Big\}\\
&\le \mathbb{E}_{t}\bigg[\exp\bigg\{2\gamma |\xi|+2\gamma |Y_{t}^{m}|+2\gamma\int_{t}^{T}\Big(\zeta_{s}(\omega)+K\big(|\Delta_{\theta}Y_{s}^{m-1,p}|+\mathbb{E}[|\Delta_{\theta}Y_{s}^{m-1,p}|]+|\Delta_{\theta}\widetilde{Y}_{s}^{m-1,p}|\\
&\hspace{2.6em}+\mathbb{E}[|\Delta_{\theta}\widetilde{Y}_{s}^{m-1,p}|]\big)+2K\big(|Y_{s}^{m+p-1}|+\mathbb{E}[|Y_{s}^{m+p-1}|]+|Y_{s}^{m-1}|+\mathbb{E}[|Y_{s}^{m-1}|]\big)\Big)ds\bigg\}\bigg].
\end{aligned}
\end{equation*}
Thus, applying Jensen's inequality once again, we observe that, for $t\in[0,T]$,
\begin{equation*}\label{eq 5.22}\tag{5.22}
\begin{aligned}
&\hspace{2em}\exp\Big\{\gamma\Big(|\Delta_{\theta}Y_{t}^{m,p}|+|\Delta_{\theta}\widetilde{Y}_{t}^{m,p}|\Big)\Big\}\\
&\le \mathbb{E}_{t}\bigg[\exp\bigg\{4n\gamma \big(|\xi|+ |Y_{t}^{m}|+|Y_{t}^{m+p}|\big)+4n\gamma\int_{t}^{T}\Big(\zeta_{s}(\omega)+K\big(|\Delta_{\theta}Y_{s}^{m-1,p}|+\mathbb{E}[|\Delta_{\theta}Y_{s}^{m-1,p}|]\\
&\hspace{0.8em}+|\Delta_{\theta}\widetilde{Y}_{s}^{m-1,p}|+\mathbb{E}[|\Delta_{\theta}\widetilde{Y}_{s}^{m-1,p}|]\big)+2K\big(|Y_{s}^{m+p-1}|+\mathbb{E}[|Y_{s}^{m+p-1}|]+|Y_{s}^{m-1}|+\mathbb{E}[|Y_{s}^{m-1}|]\big)\Big)ds\bigg\}\bigg].
\end{aligned}
\end{equation*}
In view of (\ref{eq 5.22}), applying Doob's maximal inequality, as well as H\"{o}lder's and Jensen's inequalities, we get that, for every $q>1$ and $t\in[0,T]$,
\begin{equation*}\label{eq 5.23}\tag{5.23}
\begin{aligned}
&\hspace{2em}\mathbb{E}\bigg[\exp\bigg\{q\gamma\sup\limits_{s\in[t,T]}\Big(|\Delta_{\theta}Y_{s}^{m,p}|+|\Delta_{\theta}\widetilde{Y}_{s}^{m,p}|\Big)\bigg\}\bigg]\\
&\le \sqrt{\widetilde{I}(q)}\bigg(\mathbb{E}\bigg[\exp\bigg\{16nq\gamma K\sup\limits_{s\in[t,T]}\big(|\Delta_{\theta}Y_{s}^{m-1,p}|+|\Delta_{\theta}\widetilde{Y}_{s}^{m-1,p}|\big)(T-t)\bigg\}\bigg]\bigg)^{\frac{1}{2}},
\end{aligned}
\end{equation*}
where
\begin{equation*}
\begin{aligned}
&\widetilde{I}(q):=R(q)\sup\limits_{m,p\ge1}\mathbb{E}\bigg[\exp\bigg\{8nq\gamma \Big(|\xi|+\sup\limits_{t\in[0,T]}|Y_{t}^{m}|+\sup\limits_{t\in[0,T]}|Y_{t}^{m+p}|\Big)\\
&\hspace{6em}+8nq\gamma\int_{0}^{T}\Big(\zeta_{s}(\omega)+2K\big(|Y_{s}^{m+p-1}|+\mathbb{E}[|Y_{s}^{m+p-1}|]+|Y_{s}^{m-1}|+\mathbb{E}[|Y_{s}^{m-1}|]\big)\Big)ds\bigg\}\bigg]
\end{aligned}
\end{equation*}
with $R(q)$ defined as specified in (\ref{eq 5.6}).

\noindent\textbf{Step 5.}\ \ We prove that, for any $q>1$, $\{(Y^{m},Z^{m})\}_{m=1}^{\infty}$ is a Cauchy sequence in $\mathcal{S}^{q}(0,T;\mathbb{R}^{n})\times \mathcal{H}^{q}(0,T;\mathbb{R}^{n\times d})$. By comparing (\ref{eq 5.9}) and (\ref{eq 5.23}), we see that, using an iteration w.r.t. $m$, a similar argument as that for (\ref{eq 5.5}) shows the existence of a positive constant $I^{\ast}(q)$ which depends only on $q$ such that, for every $q>1$,
\begin{equation*}
\begin{aligned}
&\hspace{1.8em}\mathbb{E}\bigg[\exp\bigg\{q\gamma\sup\limits_{t\in[0,T]}\Big(|\Delta_{\theta}Y_{t}^{m,p}|+|\Delta_{\theta}\widetilde{Y}_{t}^{m,p}|\Big)\bigg\}\bigg]\\
&\le I^{\ast}(q)\bigg(\mathbb{E}\bigg[\exp\bigg\{2^{[16nK T]}q\gamma\sup\limits_{t\in[0,T]}\Big(|\Delta_{\theta}Y_{t}^{1,p}|+|\Delta_{\theta}\widetilde{Y}_{t}^{1,p}|\Big)\bigg\}\bigg]\bigg)^{\frac{1}{2^{m-1}}}.
\end{aligned}
\end{equation*}
This combined with (\ref{eq 5.5}) yields that, for all $q>1$ and $\theta\in(0,1)$,
\begin{equation*}\label{eq 5.24}\tag{5.24}
\begin{aligned}
\limsup\limits_{m\rightarrow \infty}\sup\limits_{p\ge1}\mathbb{E}\bigg[\exp\bigg\{q\gamma\sup\limits_{t\in[0,T]}\Big(|\Delta_{\theta}Y_{t}^{m,p}|+|\Delta_{\theta}\widetilde{Y}_{t}^{m,p}|\Big)\bigg\}\bigg]\le I^{\ast}(q).
\end{aligned}
\end{equation*}
We now show that $\{(Y^{m},Z^{m})\}_{m=1}^{\infty}$ is a Cauchy sequence; its limit is a solution. For this end, we derive from (\ref{eq 5.24}) that, for each $\theta\in(0,1)$,
\begin{equation*}
\begin{aligned}
\limsup\limits_{m\rightarrow \infty}\sup\limits_{p\ge1}\mathbb{E}\Big[\sup\limits_{t\in[0,T]}|Y_{t}^{m+p}-\theta Y_{t}^{m}|\Big]\le (1-\theta)\frac{I^{\ast}(2)}{2\gamma}.
\end{aligned}
\end{equation*}
Thus, we have
\begin{equation*}
\begin{aligned}
\limsup\limits_{m\rightarrow \infty}\sup\limits_{p\ge1}\mathbb{E}\Big[\sup\limits_{t\in[0,T]}|Y_{t}^{m+p}-Y_{t}^{m}|\Big]\le (1-\theta)\bigg(\frac{I^{\ast}(2)}{2\gamma}+\sup\limits_{m\ge1}\mathbb{E}\Big[\sup\limits_{t\in[0,T]}|Y_{t}^{m}|\Big]\bigg)<+\infty.
\end{aligned}
\end{equation*}
By sending $\theta$ to $1$ and considering (\ref{eq 5.5}), we can deduce that there exists an adapted process $Y\in \mathcal{E}(0,T;\mathbb{R}^{n})$ such that, for all $q>1$,
\begin{equation}\label{eq 5.25}\tag{5.25}
\lim\limits_{m\rightarrow \infty}\mathbb{E}\bigg[\sup\limits_{t\in[0,T]}|Y_{t}^{m}-Y_{t}|^{q}\bigg]=0\ \ \text{and}\ \ \lim\limits_{m\rightarrow \infty}\mathbb{E}\bigg[\exp\bigg\{q\sup\limits_{t\in[0,T]}|Y_{t}^{m}-Y_{t}|\bigg\}\bigg]=1.
\end{equation}
Moreover, by applying It\^{o}'s formula to $|Y_{t}^{m+p}-Y_{t}^{m}|^{2}$, it can be concluded that, for all $m,p\ge1$,
\begin{equation}\label{eq 5.26}\tag{5.26}
\begin{aligned}
&\mathbb{E}\bigg[\int_{0}^{T}|Z_{s}^{m+p}-Z_{s}^{m}|^{2}ds\bigg]\le 2\mathbb{E}\bigg[\sup\limits_{t\in[0,T]}|Y_{t}^{m+p}-Y_{t}^{m}|\\
&\hspace{4em}\cdot\int_{0}^{T}\sum_{i=1}^{n}\Big|f^{i}(s,Y_{s}^{m+p-1},Z_{s}^{m+p;i},\mathbb{P}_{Y_{s}^{m+p-1}})-f^{i}(s,Y_{s}^{m-1},Z_{s}^{m;i},\mathbb{P}_{Y_{s}^{m-1}})\Big|ds\bigg].
\end{aligned}
\end{equation}
By using assumption (A8), as well as (\ref{eq 5.5}) and (\ref{eq 5.16}), we get that
\begin{equation}\label{eq 5.27}\tag{5.27}
\sup\limits_{m,p\ge1}\mathbb{E}\bigg[\bigg(\int_{0}^{T}\sum_{i=1}^{n}\Big|f^{i}(s,Y_{s}^{m+p-1},Z_{s}^{m+p;i},\mathbb{P}_{Y_{s}^{m+p-1}})-f^{i}(s,Y_{s}^{m-1},Z_{s}^{m;i},\mathbb{P}_{Y_{s}^{m-1}})\Big|ds\bigg)^{2}\bigg]<+\infty.
\end{equation}
Then, applying H\"{o}lder's inequality to (\ref{eq 5.26}) and combining (\ref{eq 5.25}) and (\ref{eq 5.27}) we obtain that
\begin{equation*}
\lim\limits_{m\rightarrow\infty}\sup_{p\ge1}\mathbb{E}\bigg[\int_{0}^{T}|Z_{s}^{m+p}-Z_{s}^{m}|^{2}ds\bigg]=0,
\end{equation*}
from which together with (\ref{eq 5.16}) it can be deduced that there exists $Z\in \mathcal{M}(0,T;\mathbb{R}^{n\times d})$ such that
\begin{equation*}\label{eq 5.28}\tag{5.28}
\lim\limits_{m\rightarrow\infty}\mathbb{E}\bigg[\bigg(\int_{0}^{T}|Z_{s}^{m}-Z_{s}|^{2}ds\bigg)^{\frac{q}{2}}\bigg]=0,\ \ q>1.
\end{equation*}
Finally, in view of (\ref{eq 5.25}) and (\ref{eq 5.28}), by letting $n\rightarrow \infty$ in BSDE (\ref{eq 5.4}), one can easily see that $(Y,Z)$ is a solution of system of diagonally quadratic mean-field BSDE (\ref{eq 5.1}).

\noindent\textbf{Step 6.}\ \ We prove the uniqueness of the solution. In order to obtain the uniqueness result, we assume that also $(\widetilde{Y},\widetilde{Z})\in \mathcal{E}(0,T;\mathbb{R}^{n})\times\mathcal{M}(0,T;\mathbb{R}^{n\times d})$ is a solution of mean-field BSDE (\ref{eq 5.1}). For $\theta\in(0,1)$, we define
\begin{equation*}
\Delta_{\theta}P:=\frac{Y-\theta\widetilde{Y}}{1-\theta},\ \ \Delta_{\theta}Q:=\frac{Z-\theta\widetilde{Z}}{1-\theta}\ \ \text{and}\ \ \Delta_{\theta}\widetilde{P}:=\frac{\widetilde{Y}-\theta Y}{1-\theta},\ \ \Delta_{\theta}\widetilde{Q}:=\frac{\widetilde{Z}-\theta Z}{1-\theta}.
\end{equation*}
By making use of a similar reasoning as in (\ref{eq 5.18}) to (\ref{eq 5.23}), one can easily show that, for all $q>1$,
\begin{equation*}\label{eq 5.29}\tag{5.29}
\begin{aligned}
&\hspace{1em}\mathbb{E}\bigg[\exp\bigg\{q\gamma\sup\limits_{s\in[t,T]}\Big(|\Delta_{\theta}P_{s}|+|\Delta_{\theta}\widetilde{P}_{s}|\Big)\bigg\}\bigg]\\
&\le \sqrt{\bar{I}(q)}\bigg(\mathbb{E}\bigg[\exp\bigg\{16nq\gamma K\sup\limits_{s\in[t,T]}\big(|\Delta_{\theta}P_{s}|+|\Delta_{\theta}\widetilde{P}_{s}|\big)(T-t)\bigg\}\bigg]\bigg)^{\frac{1}{2}},\ \ t\in[0,T],
\end{aligned}
\end{equation*}
where
\begin{equation*}
\begin{aligned}
&\bar{I}(q):=R(q)\mathbb{E}\bigg[\exp\bigg\{8nq\gamma \Big(|\xi|+\sup\limits_{t\in[0,T]}|Y_{t}|+\sup\limits_{t\in[0,T]}|\widetilde{Y}_{t}|\Big)\\
&\hspace{6em}+8nq\gamma\int_{0}^{T}\Big(\zeta_{s}(\omega)+2K\big(|Y_{s}|+\mathbb{E}[|Y_{s}|]+|\widetilde{Y}_{s}|+\mathbb{E}[|\widetilde{Y}_{s}|]\big)\Big)ds\bigg\}\bigg]
\end{aligned}
\end{equation*}
with $R(q)$ defined as specified in (\ref{eq 5.6}).

For the case of $K=0$, it can be deduced from (\ref{eq 5.29}) that
\begin{equation*}
\mathbb{E}\bigg[2\gamma\sup\limits_{t\in[0,T]}|\Delta_{\theta}P_{t}|\bigg]\le\mathbb{E}\bigg[\exp\bigg\{2\gamma\sup\limits_{t\in[0,T]}|\Delta_{\theta}P_{t}|\bigg\}\bigg]\le \sqrt{\bar{I}(2)}.
\end{equation*}
Then we have
\begin{equation*}
\mathbb{E}\bigg[\sup\limits_{t\in[0,T]}|Y_{t}-\widetilde{Y}_{t}|\bigg]\le (1-\theta)\bigg(\frac{\sqrt{\bar{I}(2)}}{2\gamma}+\mathbb{E}\bigg[\exp\bigg\{\sup\limits_{t\in[0,T]}|\widetilde{Y}_{t}|\bigg\}\bigg]\bigg),
\end{equation*}
and taking the limit $\theta\rightarrow 1$ implies that $Y=\widetilde{Y}$ and $Z=\widetilde{Z}$ on $[0,T]$. Otherwise, for the case of $T\le \varepsilon^{\ast}:=\frac{1}{16nK}$, it follows from (\ref{eq 5.29}) that
\begin{equation*}
\mathbb{E}\bigg[2\gamma\sup\limits_{t\in[0,T]}|\Delta_{\theta}P_{t}|\bigg]\le\mathbb{E}\bigg[\exp\bigg\{2\gamma\sup\limits_{t\in[0,T]}|\Delta_{\theta}P_{t}|\bigg\}\bigg]\le \bar{I}(2),
\end{equation*}
and then $Y=\widetilde{Y}$ and $Z=\widetilde{Z}$ on $[0,T]$. Similarly, let $n_{0}$ be the unique positive integer such that $16nKT \le n_{0} < 16nKT+1$. Then we can successively obtain the uniqueness on the time intervals $[T-i\varepsilon^{\ast},T-(i-1)\varepsilon^{\ast}],\ 1\le i\le n_{0}-1$ and $[0,T-(n_{0}-1)\varepsilon^{\ast}]$. Consequently, we have established the uniqueness on the whole interval $[0,T]$, and the proof of Theorem 5.1 is hereby finished.
\end{proof}

At the end of this section, we consider the following Volterra-type mean-field BSDE:
\begin{equation*}\label{eq 5.30}\tag{5.30}
Y_{t}=\xi+\int_{t}^{T}\mathbb{E}_{t}\Big[g(s,Y_{s\vee\cdot},Z,\mathbb{P}_{(Y_{s\vee \cdot},Z)})\Big]ds+\int_{t}^{T}f(s,Z_{s})ds-\int_{t}^{T}Z_{s}dW_{s},\ t\in[0,T],
\end{equation*}
where $g:\Omega\times[0,T]\times \mathcal{C}([0,T];\mathbb{R}^{n})\times L^{2}([0,T];\mathbb{R}^{n\times d})\times \mathcal{P}_{2}(\mathcal{C}([0,T];\mathbb{R}^{n})\times L^{2}([0,T];\mathbb{R}^{n\times d}))\rightarrow \mathbb{R}^{n}$ and $f:\Omega\times[0,T]\times\mathbb{R}^{n\times d}\rightarrow \mathbb{R}^{n}$ are $\mathbb{F}$-progressively measurable and satisfy the following assumptions:\vspace{1em}

\noindent\textbf{(A11)} For $1\le i\le n,\ f^{i}(\omega,t,z)$ varies with $(\omega,t,z^{i})$ only, and grows quadratically in $z^{i}$, i.e., $d\mathbb{P}\times dt\text{-}$almost all $(\omega,t)\in \Omega\times [0,T]$, there exists $\gamma>0$ and an $\mathbb{F}$-progressively measurable process $\zeta=(\zeta_{t})$ with $\mathbb{E}\big[\exp\big\{p\int_{0}^{T}\zeta_{t}(\omega)dt\big\}\big]<+\infty$, for all $p\ge1$, such that, for all $z\in\mathbb{R}^{n\times d}$,
\begin{equation*}
|f^{i}(\omega,t,z)|\le \zeta_{t}(\omega)+\frac{\gamma}{2}|z^{i}|^{2}.
\end{equation*}

\noindent\textbf{(A12)} For $1\le i\le n$, it holds that $f^{i}(\omega,t,\cdot)$ is either convex or concave.

\noindent\textbf{(A13)} There exists $C\ge0$ such that for all $t\in[0,T]$, $y,\bar{y}\in \mathcal{C}([0,T];\mathbb{R}^{n})$, $z,\bar{z}\in L^{2}([0,T];\mathbb{R}^{n\times d})$ and $\mu,\bar{\mu}\in \mathcal{P}_{2}(\mathcal{C}([0,T];\mathbb{R}^{n})\times L^{2}([0,T];\mathbb{R}^{n\times d}))$ with $\mu(\mathcal{C}([0,T];\mathbb{R}^{n})\times \cdot)=\bar{\mu}(\mathcal{C}([0,T];\mathbb{R}^{n})\times \cdot)$,
\begin{equation*}
\begin{aligned}
&|g(t,y,z,\mu)|\le C\big(1+|y|_{\mathcal{C}_{T}}+|z|_{L_{T}^{2}}+W_{2}(\mu_{1},\delta_{0})+W_{2}(\mu_{2},\delta_{0})^{2}\big),\\
&|g(t,y,z,\mu)-g(t,\bar{y},z,\bar{\mu})|\le C\big(|y-\bar{y}|_{\mathcal{C}_{T}} +W_{2}(\mu_{1},\bar{\mu}_{1})\big).\\
\end{aligned}
\end{equation*}

\begin{theorem} Let the assumptions (A7) and (A11)-(A13) be satisfied. Then mean-field BSDE (\ref{eq 5.30}) has a solution $(Y,Z)\in \mathcal{S}^{2}_{\mathbb{F}}(0,T;\mathbb{R}^{n})\times \mathcal{M}(0,T;\mathbb{R}^{n\times d})$ on the whole interval $[0,T]$. Moreover, if the function $g$ is bounded, the solution $(Y,Z)$ belongs to $\mathcal{E}(0,T;\mathbb{R}^{n})\times \mathcal{M}(0,T;\mathbb{R}^{n\times d})$, and it is unique in this space.
\end{theorem}
\begin{proof}
Firstly, by virtue of Theorem 5.1 and assumptions (A7) and (A11)-(A12), we have the following quadratic BSDE with an unbounded terminal value admits a unique solution $(Y^{\prime},Z)\in \mathcal{E}(0,T;\mathbb{R}^{n})\times \mathcal{M}(0,T;\mathbb{R}^{n\times d})$:
\begin{equation*}\label{eq 5.31}\tag{5.31}
Y_{t}^{\prime}=\xi+\int_{t}^{T}f(s,Z_{s})ds-\int_{t}^{T}Z_{s}dW_{s},\ t\in[0,T].
\end{equation*}
Furthermore, we consider the following equation:
\begin{equation*}\label{eq 5.32}\tag{5.32}
Y_{t}=Y_{t}^{\prime}+\int_{t}^{T}\mathbb{E}_{t}\Big[g(s,Y_{s\vee\cdot},Z,\mathbb{P}_{(Y_{s\vee \cdot},Z)})\Big]ds,\ t\in[0,T].
\end{equation*}
To obtain the existence of the solution, we now apply Picard iteration. Let $Y_{t}^{0}:=0,\ t\in[0,T]$, and for $n\ge 0$, we define
\begin{equation*}\label{eq 5.33}\tag{5.33}
Y_{t}^{n+1}=Y_{t}^{\prime}+\int_{t}^{T}\mathbb{E}_{t}\Big[g(s,Y_{s\vee\cdot}^{n},Z,\mathbb{P}_{(Y_{s\vee \cdot}^{n},Z)})\Big]ds,\ t\in[0,T].
\end{equation*}
By iteration we see that $Y^{n+1}\in \mathcal{S}_{\mathbb{F}}^{2}(0,T;\mathbb{R}^{n}),\ n\ge0$.

Let us put $\Delta Y^{n}:=Y^{n}-Y^{n-1}$. Then we have
\begin{equation*}\label{eq 5.34}\tag{5.34}
\Delta Y_{t}^{n+1}=\int_{t}^{T}\mathbb{E}_{t}\Big[g(s,Y_{s\vee\cdot}^{n},Z,\mathbb{P}_{(Y_{s\vee \cdot}^{n},Z)})-g(s,Y_{s\vee\cdot}^{n-1},Z,\mathbb{P}_{(Y_{s\vee \cdot}^{n-1},Z)})\Big]ds.
\end{equation*}

\noindent Applying It\^{o}'s formula to $e^{\beta t}|\Delta Y_{t}^{n+1}|^{2}$ ($\beta>0$ will be specified later) and combining it with (A13), Young's and Jensen's inequalities, we obtain
\begin{equation*}\label{eq 5.35}\tag{5.35}
\begin{aligned}
&\hspace{2.6em}e^{\beta t}|\Delta Y_{t}^{n+1}|^{2}+\beta\int_{t}^{T}e^{\beta s}|\Delta Y_{s}^{n+1}|^{2}ds\\
&\le\frac{\beta}{2}\int_{t}^{T}e^{\beta s}|\Delta Y_{s}^{n+1}|^{2}ds+\frac{4C^{2}}{\beta}\int_{t}^{T}\Big(\mathbb{E}_{t}[\sup_{r\in[s,T]}e^{\beta r}|\Delta Y_{r}^{n}|^{2}]+\mathbb{E}[\sup_{r\in[s,T]}e^{\beta r}|\Delta Y_{r}^{n}|^{2}]\Big)ds.
\end{aligned}
\end{equation*}

\noindent Taking the supremum over $t\in[0,T]$ in (\ref{eq 5.35}) and then the expectation, combined with Fubini's theorem, this gives
\begin{equation*}
\mathbb{E}[\sup\limits_{t\in[0,T]}e^{\beta t}|\Delta Y_{t}^{n+1}|^{2}]\le \frac{8C^{2}T}{\beta}\mathbb{E}[\sup\limits_{t\in[0,T]}e^{\beta t}|\Delta Y_{t}^{n}|^{2}].
\end{equation*}

\noindent By choosing $\beta=32C^{2}T$, we obtain
\begin{equation*}
\mathbb{E}[\sup\limits_{t\in[0,T]}e^{\beta t}|\Delta Y_{t}^{n+1}|^{2}]\le \frac{1}{4}\mathbb{E}[\sup\limits_{t\in[0,T]}e^{\beta t}|\Delta Y_{t}^{n}|^{2}].
\end{equation*}

\noindent By induction, we can derive that there is a positive constant $M$ such that, for all $n\ge 1$,
\begin{equation*}
\mathbb{E}[\sup\limits_{t\in[0,T]}|\Delta Y_{t}^{n}|^{2}]\le \mathbb{E}[\sup\limits_{t\in[0,T]}e^{\beta t}|\Delta Y_{t}^{n}|^{2}]\le \frac{1}{4^{n-1}}\mathbb{E}[\sup\limits_{t\in[0,T]}e^{\beta t}|\Delta Y_{t}^{1}|^{2}]\le \frac{M^{2}}{4^{n}}.
\end{equation*}

\noindent Thus, for $m>n$,
\begin{equation*}
\Vert Y^{m}-Y^{n}\Vert_{\mathcal{S}_{\mathbb{F}}^{2}(0,T)}\le \sum\limits_{k=n+1}^{m}\Vert\Delta Y^{k}\Vert_{\mathcal{S}_{\mathbb{F}}^{2}(0,T)}\le \sum\limits_{k=n+1}^{m}\frac{M}{2^{k}}\le \frac{M}{2^{n}}\rightarrow 0,\ \ \text{as}\ n\rightarrow +\infty.
\end{equation*}

\noindent Therefore, we conclude that there exists $Y\in\mathcal{S}_{\mathbb{F}}^{2}(0,T;\mathbb{R}^{n}) $ such that
\begin{equation*}
\mathbb{E}\Big[\sup\limits_{t\in[0,T]}|Y_{t}^{n}-Y_{t}|^{2}\Big]\rightarrow 0,\ \ \text{as}\ n\rightarrow +\infty.
\end{equation*}
Now by taking the limit as $n\rightarrow +\infty$ in (\ref{eq 5.33}), it can be easily confirmed that $Y$ satisfies (\ref{eq 5.32}).

In conclusion, from (\ref{eq 5.31})-(\ref{eq 5.32}) we can get that $(Y,Z)\in \mathcal{S}_{\mathbb{F}}^{2}(0,T;\mathbb{R}^{n})\times \mathcal{M}(0,T;\mathbb{R}^{n\times d})$ is a solution of quadratic mean-field BSDE (\ref{eq 5.30}).

Moreover, if $g$ is bounded, the fact that $(Y^{\prime},Z)\in\mathcal{E}(0,T;\mathbb{R}^{n})\times \mathcal{M}(0,T;\mathbb{R}^{n\times d})$ implies that also $(Y,Z)$ belongs to this space. Again under the assumption that $g$ is bounded, if $(\widetilde{Y},\widetilde{Z})\in\mathcal{E}(0,T;\mathbb{R}^{n})\times \mathcal{M}(0,T;\mathbb{R}^{n\times d})$ is a solution of (\ref{eq 5.30}), also $(\widetilde{Y}^{\prime}=(\widetilde{Y}^{\prime}_{t}:=\widetilde{Y}_{t}-\int_{t}^{T}\mathbb{E}_{t}\big[g(s,\widetilde{Y}_{s\vee\cdot},\widetilde{Z},\mathbb{P}_{(\widetilde{Y}_{s\vee \cdot},\widetilde{Z})})\big]ds,\ t\in[0,T]),\widetilde{Z})$ belongs to $\mathcal{E}(0,T;\mathbb{R}^{n})\times \mathcal{M}(0,T;\mathbb{R}^{n\times d})$ and solves BSDE (\ref{eq 5.31}). However, (\ref{eq 5.31}) has a unique solution in $\mathcal{E}(0,T;\mathbb{R}^{n})\times \mathcal{M}(0,T;\mathbb{R}^{n\times d})$. Hence, $(\widetilde{Y}^{\prime},\widetilde{Z})=(Y^{\prime},Z)$. The fact that $\widetilde{Y}=Y$ follows now from the uniqueness of equation (\ref{eq 5.32}).
\end{proof}

\begin{remark}
Theorem 5.2 still holds true for the following special form of mean-field BSDE (\ref{eq 5.30}):
\begin{equation*}
Y_{t}=\xi+\int_{t}^{T}g(s,Y_{s},Z_{s},\mathbb{P}_{(Y_{s},Z_{s})})ds+\int_{t}^{T}f(s,Z_{s})ds-\int_{t}^{T}Z_{s}dW_{s},\ t\in[0,T].
\end{equation*}
\end{remark}

\end{document}